\documentclass[a4paper,10pt]{article}

\usepackage[
a4paper,
top=2cm,
bottom=2cm,
left=2cm,
right=2cm
]{geometry}

\usepackage[numbers]{natbib}
\usepackage{amsmath,amsthm}
\usepackage{amssymb}
\usepackage{mathrsfs}
\usepackage{mathpazo}
\usepackage{graphicx}
\usepackage{enumerate}
\usepackage[utf8,latin1]{inputenc}
\usepackage[T1]{fontenc} 

\usepackage{mathtools}
\usepackage{xcolor}
\usepackage{caption}
\usepackage{subcaption}
\usepackage{tabularx}
\usepackage{multirow}
\usepackage{booktabs}
\usepackage{array}
\usepackage[percent]{overpic}
\usepackage{bm}
\usepackage{algorithm}
\usepackage{algpseudocode}
\usepackage{algorithm}
\usepackage{siunitx}
\usepackage{amsopn}
\usepackage{etoolbox}
\usepackage{hyperref}

\newtheorem{theorem}{Theorem}
\newtheorem{definition}[theorem]{Definition}
\newtheorem{corollary}[theorem]{Corollary}
\newtheorem{rmk}{Remark}
\newcommand{\rom}[1]{%
	\textup{\uppercase\expandafter{\romannumeral#1}}%
}

\DeclareMathOperator{\trace}{trace}

\DeclareMathOperator{\diag}{diag}
\DeclareMathOperator{\Fourier}{{\mathscr F}}

\title{Reconstruction Formulae for 3D Field-Free Line Magnetic Particle Imaging}

\usepackage{authblk}

\author[1,2]{Vladyslav Gapyak \thanks{Emails: vladyslav.gapyak@h-da.de, thomas.maerz@h-da.de, andreas.weinmann@h-da.de}}
\author[1,2]{Thomas M{\"a}rz}
\author[1,2]{Andreas Weinmann}
\affil[1]{Algorithms for Computer Vision, Imaging and Data Analysis Lab at Darmstadt University of Applied Sciences, Sch\"{o}fferstr. 3, 64295, Darmastdt, Germany}
\affil[2]{Data Science Institute, European University of Technology, European Union}

\date{}
\begin{document}
\maketitle

\begin{abstract}
    Magnetic Particle Imaging (MPI) is a promising noninvasive in vivo imaging modality that makes it possible to map the spatial distribution of superparamagnetic nanoparticles by exposing them to dynamic magnetic fields. In the Field-Free Line (FFL) scanner topology, the spatial encoding of the particle distribution is performed by applying magnetic fields vanishing on straight lines. The voltage induced in the receiving coils by the particles when exposed to the magnetic fields constitute the signal from which the particle distribution is to be reconstructed. To avoid lengthy calibration, model-based reconstruction formulae have been developed for the 2D FFL scanning topology. In this work we develop reconstruction formulae for 3D FFL. Moreover, we provide a model-based reconstruction algorithm for 3D FFL and we validate it with a numerical experiment.
\end{abstract}

\noindent{\it Keywords}: magnetic particle imaging, model-based reconstruction, reconstruction formulae, phase space, inverse problems, variational regularization, field-free line, Radon transform

\section{Introduction}

Magnetic Particle Imaging (MPI) is a promising medical imaging modality that allows to reconstruct the distribution of superparamagnetic nanoparticles by exploiting their non-linear magnetization response to dynamic magnetic fields. The principles behind MPI have been introduced by Gleich and Weizenecker in their seminal paper published in Nature in 2005~\cite{GleichWeizenecker2005} and led them and their team to win the European Inventor Award in 2016. In 2009, Weizenecker \emph{et al.}~\cite{Weizenecker_etal2009} successfully managed to produce a video (3D+time) of a distribution of a clinically approved MRI contrast agent flowing through the beating heart of a mouse, thereby taking an important step towards medical applications of MPI and in particular towards its applicability to 3D real-time \emph{in vivo} imaging. Since then, interest in the technology has been on the rise, as has the number of medical applications. These include multimodal imaging~\cite{aramiMultimodalImagingMPI}, cancer detection (of as few as 250 cancer cells~\cite{Song2018}) and cancer imaging~\cite{aramiMultimodalImagingMPI,Du2019,Tay2021,Yu2017}, (stem) cell tracing~\cite{connell2015advancedcellTherapies,GoodwillConolly2011,JUNG2018139,Lemaster2018,Tomitaka2015Lactoferrin}, inflammation tracing and lung perfusion imaging~\cite{Zhou_2017}, drug delivery and monitoring~\cite{Zhu2019}, cardiovascular~\cite{Bakenecker2018MPIvascular,Tong2021,Vaalma2017} and blood flow~\cite{Franke2020BloodFlow} imaging, tracking of medical instruments~\cite{haegele2012instrumentvisualization}, brain injury detection~\cite{Orendorff2017} as well as future MPI-assisted stroke detection and brain monitoring made possible by a recently developed human-sized MPI scanner~\cite{Graeser2019humanbrain}. The high degree of applicability of MPI is due to its benefits, including quatifyability and high sensitivity, and due to its advantages over other imaging modality such as CT~\cite{Buzug2008CTbook}, MRI~\cite{schamelsafety2015}, PET~\cite{PET} and SPECT~\cite{SPECT}, for example, its ability to provide higher spatial resolution in a shorter acquisition time and the absence of radiation or radioactive tracers, making it a safer imaging option~\cite{BuzugKnopp2012}. Deeper comparison between existing imaging modalities and MPI as well as further examples of state-of-the-art applications can be found in~\cite{billingsMPIapplications,Yang2022}.

In brief, the principle and the procedure behind MPI are as follows: (i) the specimen to be scanned is injected with a tracer containing the magnetic nanoparticles (MNs); (ii) a dynamic magnetic field with a low-field volume (LFV) is applied and the LFV is spatially moved; (iii) because the response of the MNs is nonlinear precisely where the field is close to zero,  the LFV acts as a sensitive region whose motion causes the particles to induce a voltage in the receiver coils, a voltage that constitutes the signal that encodes the information about the particles' position and concentration. The imaging task of MPI is the reconstruction of the spatial distribution of the particles in the tracer from the signal (voltage) acquired during a scan.

Concerning the LFVs employed, the first possibility explored~\cite{GleichWeizenecker2005} is the case when the dynamic magnetic field applied is the superposition of a static field with a point where the field vanishes - the Field-Free Point (FFP) - and a dynamic field that steers the FFP along some trajectory in the scanning region - or Field of View (FoV). There are currently two main approaches to reconstructing distributions with the FFP topology: a measurement-based approach and a model-based approach. In the measurement-based approach a 3D grid (tessellation with cubic cells) of the volume to be scanned inside the scanner's bore is considered; a probe with a reference concentration of tracer is iteratively positioned and scanned at each voxel, i.e., at each cubic cell constituting the grid. This way the response of the scanning system to discrete Dirac delta impulses (the point spread function~\cite{bertero2021introduction}) at each voxel position is collected in a system matrix, which should ideally describe the linear relationship between discrete (according to the chosen grid) distributions and the acquired signal. The reconstruction idea is, given the MPI signal collected, to invert the system matrix (with regularization techniques due to the presence of noise and the ill-conditioning of the matrix) to obtain a discretized reconstruction of the target distribution~\cite{Knopp_etal2010ec,storath2016edge}. The acquisition of the system matrix (calibration) is time-consuming and while effort has been put into the research of ways to speed up this procedure~\cite{knopp2011prediction,Lampe_etal2012,Rahmer_etal2012}, other solutions have been investigated that avoid the lengthy calibration and are model based. Indeed, model-based approaches are an active area of research in MPI and try to provide reconstruction formulae for various MPI setups. In the particular case of the FFP topology, at first direct reconstruction formulae were available only for 1D scans~\cite{Rahemeretal2009} and employed Chebyshev polynomials but were followed soon after by the so-called $X$-space formulation for the 1D~\cite{GoodwillConolly2010} scans and later for the multidimensional case with Cartesian scanning trajectories~\cite{GoodwillConolly2011}. Following this, direct reconstruction with Chebyshev polynomials has also been developed for the 2D and 3D case for Lissajous trajectories~\cite{droigk2022multidimcheb}. Using the Langevin theory, the authors derived reconstruction formulae for the 2D and 3D FFP case~\cite{marz2016model} and formulated a two-stage algorithm that has been recently further developed~\cite{gapyak2022mdpi}. A remarkable point of the reconstruction formulae in~\cite{marz2016model,gapyak2022mdpi} is that the reconstruction formulae are independent from the choice of the scanning trajectory; this property has been leveraged to perform reconstructions in multi-patch scenarios~\cite{gapyak2023multipatch} and is inherited by the methods proposed in this paper. Multi-patching is a methodology in MPI which aims at solving an MPI specific difficulty: in order to avoid high specific absorption rate (SAR) and peripheral nerves stimulation (PNS)~\cite{schamelsafety2015,Saritas2013-jn}, the size of the field of view FoV does not exceed a couple of centimeters; in order to cover bigger regions of interest with a (comparably) small FoV, multi-patching is employed by either physically moving the specimen in the scanner~\cite{Szwargulskimovingtable2018} or by considering an additional focus field~\cite{Weizenecker2010FastMD}. In 2014, Vogel et al.~\cite{vogel2014TWMPI} proposed the Traveling Wave MPI (TWMPI) approach as an alternative approach to cover a large FoV, which combines the advantages of the drive-coil~\cite{gleich2005original} setup with the focus-field approach~\cite{Weizenecker2010FastMD}.

Concerning other scanning topologies, in 2008 Weizenecker, Gleich and Borgert proposed a Field-Free Line (FFL) scan modality, i.e., the employment of magnetic fields vanishing on a line, with the promise of a more effective encoding scheme, an increase in sensitivity of acquisition, a greater signal to noise ratio (SNR)~\cite{Weizenecker_2008_FFL} and potentially a lower level of ill-posedness for the inverse problem underlying the reconstruction task in MPI with FFL as opposed to the FFP topology~\cite{Kluth_2018Illposedness}. At first, there were doubts about the feasibility of a FFL scanner, but Knopp \emph{et al.}~\cite{ErbeFFLbook,Knopp2010EfficientGO,Knopp2010FFLformation} showed that it is possible to produce a system with power dissipation comparable to that of a FFP scanner. Soon after, the first 2D images with a FFL scanner were made possible~\cite{Goodwill2012Proj}. In 2014 the first continuously rotated FFL scanner has been presented by Weber \emph{et al.}~\cite{weber2014rotFFL}. Since then, effort has been put into the study of possible human-sized FFL scanners~\cite{bringout_phdthesis} and into the design of animal-sized scanners~\cite{Choi2020ScRep,MattinglyFFL2022}. As for the FFP topology, reconstruction with FFL data can be done by either acquiring a system matrix or using model-based approaches~\cite{Ilbey2017ComparisonOS}. The measurement-based approach works in the same way as with the FFP topology: a reference concentration that plays the role of a Dirac delta distribution is placed and scanned in each voxel position and the data thus acquired form the system matrix; regularized inversion of the system matrix produces a reconstruction of the target distribution (examples of measurement-based reconstructions for 2D FFL can be found in~\cite{Ilbey2017ComparisonOS,Kilic2023,Top2019trajectoryFFL}). Much more attention have received model-based approaches for FFL, thanks to the link with the Radon transform and consequent availability of a Fourier slice theorem for model-based 2D reconstructions~\cite{Knopp_2011_Fourierslicethm}. This is beneficial because the theory of the Radon transform~\cite{helgason2011,RadonOriginal1917,markoe_2006} and its inversion in 2D is well studied in the context of CT~\cite{NattererCTbook} and offers a plethora of well-established classical techniques~\cite{Willemink2019} as well as more recent Machine Learning-based methods~\cite{Szczykutowicz2022ARO} that become available for model-based reconstructions in MPI. Thanks to recent developments of the scanner, it is possible to both rotate the FFL and translate it in 3D~\cite{Top2020FFLopensidedscanner} and this flexibility has been employed for the first simulation of FFL 3D tomographic scans. Indeed, by scanning in a slice-wise fashion~\cite{soydan20213Dslicewise} the Radon-based reconstruction formulae for 2D FFL applied to each slice can be used to reconstruct 3D distributions. Scanning trajectories for tomographic 3D FFL have also been investigated and utilized with the system-matrix-based approach~\cite{Soydan20233DtrajectoryFFL} but to the best of our knowledge, no fully 3D FFL reconstruction formulae exists to date. Interestingly, the TWMPI approach can be applied to the FFL scanning topology as well~\cite{greiner2022twmpiFFL} and lead Vogel et al. to the development in 2023 of a portable human-sized scanner, the iMPI~\cite{vogel2023impi}. The iMPI scanner is, to our view, a candidate to implement the scanning set-up proposed and it uses system-matrix-based reconstruction (either obtained with a calibration procedure or simulated with model-based principles). Before describing the contributions of this paper and its outline, we would like to briefly discuss a further reason why model-based approaches for FFL scans in 3D setups are an interesting and important direction of research. In all generality, the Radon transform is a transform that considers integrals on hyperplanes in spaces of any dimension~\cite{markoe_2006} and reconstruction formulae analogous to the 2D case exist for the 3D case if we consider integrals on planes~\cite{NattererCTbook}. In MPI this translates into considering magnetic fields with a Field-Free Plane, instead of a FFL. However, a Field-Free Plane is not possible due to Maxwell's equations (a proof of this fact can be found in Section~\ref{sec:models}). As a consequence, the maximum dimensionality of possible FFRs (Field-Free Regions) in MPI is 1, which corresponds to the FFL topology and this fact underlines the importance of model-based approaches for 3D FFL scans.

\paragraph{Contributions} The purpose of this paper is to provide model-based reconstruction approaches for a proposed 3D FFL scanning scheme. In particular, we make the following major contributions:
\begin{itemize}
	\item[(i)] we provide model-based reconstruction formulae for 3D FFL MPI;
	\item[(ii)] we provide a reconstruction algorithm for 3D FFL and illustrate its applicability with a numerical example.
\end{itemize}

\paragraph{Outline of the Paper}

We begin our work by recalling the mathematical modeling of MPI in Section \ref{sec:models}, where we also describe the possible scanner topologies and prove that Field-Free Plane topologies are not possible as a consequence of Maxwell's equations. In Section \ref{sec:3d:ffl} we describe the proposed 3D FFL scanning setup and we prove the main result of this paper from which we derive the reconstruction formulae. In Section \ref{sec:rec:algorithm} we describe the three steps of the proposed reconstruction algorithm and in Section \ref{sec:discussion} we discuss the feasibility of the proposed scenario with modified real scanners, in particular with the iMPI scanner \cite{vogel2023impi}. Finally, in Section \ref{sec:experiments} we show a numerical experiment that illustrates the reconstruction algorithm. We conclude with a discussion of the results and future research directions in Section \ref{sec:conclusions}.

\section{Particle Magnetization, Induced Signal and Scanner Topologies}\label{sec:models}

We now describe the physical model that puts the magnetization response of the particles when exposed to dynamic magnetic fields in relation to the voltage induced in the receiving coils of the scanner.

\paragraph{Particle Magnetization}

Let $\rho \colon\mathbb{R}^3\to\mathbb{R}^+$ be a distribution of superparamagnetic particles exposed to a dynamic magnetic field $\bm{H}\colon\mathbb{R}^3\times\mathbb{R}\to\mathbb{R}^3$ with $(\bm{x},t)\mapsto\bm{H}(\bm{x},t)$. Assuming instantaneous remagnetization of the particles, the magnetization is modeled according to the Langevin theory of paramagnetism~\cite{chikazumi1978physics,jiles1998introduction} as follows:
\begin{equation}\label{eq:magnetization}
	\bm{M}(\bm{x},t) = m\rho (\bm{x})\mathcal{L}\biggl (\frac{\lVert \bm{H}(\bm{x},t)\rVert}{H_{\mathrm{sat}}}\biggr )\frac{\bm{H}(\bm{x},t)}{\lVert \bm{H}(\bm{x},t)\rVert} ,\quad\text{with }H_{\mathrm{sat}} =\frac{k_B T}{\mu_0 M_{\mathrm{sat}}\frac{\pi}{6}d^3},
\end{equation}
where $m$ is the magnetic moment of a single particle, $k_B$ is the Boltzmann constant, $T$ the temperature of the particles, $\mu_0$ is the magnetic permeability, $M_{\mathrm{sat}}$ the saturation magnetization, $d$ is the particles' diameter and $\mathcal{L}$ is the Langevin function $\mathcal{L}(x) = \coth (x)-\frac{1}{x}$.

\paragraph{Signal Encoding}

Let now $\bm{u}(t)$ be the voltage received by the three receive coils and $P\colon\mathbb{R}^{3}\to\mathbb{R}^{3\times 3}$ be their sensitivity pattern. According to Faraday's law of induction, the voltage induced in the coils is the negative rate of change of the magnetic flux $\bm{\Phi}(t)$:
\begin{equation}\label{eq:faraday}
	\bm{u}(t)=-\frac{d}{dt}\bm{\Phi}(t),\quad\text{where }\bm{\Phi}(t)=\mu_0 \int_{\mathbb{R}^{3}} P(\bm{x}) \bigl (\bm{H}(\bm{x},t)+\bm{M}(\bm{x},t)\bigr )\, d\bm{x} .
\end{equation}
Because only the magnetization term $\bm{M}$ in Equation \eqref{eq:faraday} is related to the distribution $\rho$, see Equation \eqref{eq:magnetization}, after subtracting from $\bm{u}(t)$ the signal of an empty scan, i.e., a scan with $\rho (\bm{x})= 0$, we obtain the signal
\begin{equation}\label{eq:signal:encoding}
	\bm{s}(t)=-\mu_0 \frac{d}{dt}\int_{\mathbb{R}^3}P(\bm{x})\bm{M}(\bm{x},t)\, d\bm{x}.
\end{equation}

\paragraph{Possible Scanner Topologies}

From \eqref{eq:signal:encoding} and \eqref{eq:magnetization} it is clear that the signal $\bm{s}(t)$ induced in the receiving coils during a scan depends on the derivative of the Langevin function, which is a bell-shaped function with maximum in $0$ and asymptotically zero away from $0$. This means that the regions that produce the highest amount of signal are those close to the set of points where the magnetic field vanishes, the so-called Field-Free Region (FFR). The shape of the FFR defines the topology (or geometry) of the MPI scanner employed. The MPI principle has been introduced by Weizenecker and Gleich~\cite{GleichWeizenecker2005} considering magnetic fields vanishing in a single point. This is known as the Field-Free Point (FFP) topology and model-based reconstruction formulae in 2D and 3D have been proposed by the authors~\cite{marz2016model}.

It is also possible to consider the Field-Free Line (FFL) topology, i.e., scanners employing magnetic fields vanishing along a line. The FFL topology has been introduced in~\cite{Weizenecker_2008_FFL} and has been tightly connected to the theory of the Radon transform in 2D and in particular, a Fourier slice theorem has been proved and provides reconstruction formulas for 2D scans~\cite{Knopp_2011_Fourierslicethm}. The Radon transform~\cite{RadonOriginal1917,NattererCTbook} is an integral transform that represents functions $f\colon\mathbb{R}^n\to\mathbb{R}$ (in a suitable space) through the integrals of $f$ on affine hyperplanes $\bm{\Pi}$ of $\mathbb{R}^n$,i.e., on affine subspaces of dimension $n-1$. One then would expect that a natural generalization of the 2D Field Free Line topology in 3D would be a Field Free Plane topology. The authors in~\cite{Weizenecker_2008_FFL} have mentioned that it is not possible to have magnetic fields vanishing on a plane, as Maxwell's equations~\cite{jackson_classical_1999} put physical limitations to the magnetic fields that can occur. We have not however found a proof of this statement in the literature and we provide one to formally motivate our choice to focus on the FFL topology and set aside the study of the Radon transform on planes in the context of MPI. The result is based on the following property of sinusoidal and irrotational fields:
\begin{theorem}\label{thm:no:ffplanes}
	Let $\Omega\subset\mathbb{R}^3$ be a non empty and connected open set, $\bm{H}\colon \Omega\to\mathbb{R}^3$ be a continuously differentiable vector field satisfying the differential equations:
	\begin{align}
		\nabla\cdot \bm{H}&=0 , \label{eq:zero:div} \\
		\nabla\times\bm{H}&=0 . \label{eq:zero:rot}
	\end{align}
	Then, level sets of $\bm{H}$ cannot be ($2$-dimensional) regular surfaces in $\Omega$.
\end{theorem}

\begin{proof}
	Suppose for the sake of contradiction that $\bm{H}$ is constant on some $2$-dimensional regular surface in $\Omega$, i.e., there exists a $2$-dimensional regular surface $\mathcal{S}\subset\Omega$ and a constant $\bm{c}\in\mathbb{R}^3$ such that in each point $\bm{x}_0\in\mathcal{S}$ it holds $\bm{H}(\bm{x}_0 )=\bm{c}$ and the tangent space $T_{\bm{x}_0}\mathcal{S}$ is 2-dimensional.
	Given a point $\bm{x}_0\in\mathcal{S}$, let $\varphi (\tau ,\sigma )$ be a surface parametrization of $\mathcal{S}$ with:
	\begin{equation}
		\begin{cases}
			\varphi (0,0) =\bm{x}_0 \\
			\partial_{\tau}\varphi (0,0) = \bm{v} \\
			\partial_{\sigma}\varphi (0,0) = \bm{u}
		\end{cases},
	\end{equation}
	for some linearly independent vectors $\bm{u}$ and $\bm{v}$. Let $\bm{J_H}$ be the Jacobian matrix of $\bm{H}$, $(\bm{J_H})_{i,j} = \frac{\partial H_i}{\partial x_j}$. Because $\bm{H}\left ( \varphi (\tau ,\sigma )\right ) =\bm{c}$, we get that $\bm{J_H}(\bm{x}_0 )\bm{v}= 0$ and $\bm{J_H}(\bm{x}_0 ) \bm{u}= 0$, which in turn means that $\bm{J_H}(\bm{x}_0 )$ has eigenvalue $0$ with geometric multiplicity of $2$. In turn, also the algebraic multiplicity must be at least $2$, meaning that if  $\lambda_1$, $\lambda_2$, $\lambda_3$ are the eigenvalues of $\bm{J_H}(\bm{x}_0 )$, then two of them are zero. Without loss of generality, we can consider $\lambda_1 = \lambda_2 = 0$. Now, we combine the hypothesis in \eqref{eq:zero:div} with the fact that the divergence of the field is the trace of the Jacobian, i.e.,
	\begin{equation}
		0 = \nabla\cdot\bm{H}(\bm{x}_0 )= \trace\left ( \bm{J_H}(\bm{x}_0 )\right ) = \lambda_1 + \lambda_2 + \lambda_3 = 0 + 0 +\lambda_3 ,
	\end{equation}
	which forces $\lambda_3 = 0$ as well and shows that $0$ is an eigenvalue with algebraic multiplicity of $3$.
	
	On the other hand, \eqref{eq:zero:rot} implies that $\bm{J_H}(\bm{x}_0 )$ is symmetric and therefore diagonalizable. Diagonalizability of a matrix is equivalent to the fact that the algebraic and geometric multiplicities of each eigenvalue are equal, which proves that the geometric multiplicity of $\bm{J_H}(\bm{x}_0 )$ is also of $3$ and consequently, that $\bm{J_H}(\bm{x}_0 )= 0$. Finally, on a regular surface the tangent space \footnote{The tangent space $T_{\bm{x}_0}\mathcal{S}$ to $\mathcal{S}$ in a point $\bm{x}_0$ is the set of velocities of paths through $\bm{x}_0$ on $\mathcal{S}$, i.e., $T_{\bm{x}_0}\mathcal{S}=\lbrace \gamma '(0) | \gamma\colon (-\varepsilon ,\varepsilon)\to\mathcal{S}$ s.t. $\gamma (0)=\bm{x}_0\rbrace$.} \cite{docarmo2016crves_surfaces} is isomorphic to the kernel of the Jacobian, i.e., $\mathrm{ker}\bm{J_H}(\bm{x}_0 )\cong T_{\bm{x}_0}\mathcal{S}$, but this leads to a contradiction: the tangent space on a $2$-dimensional regular surface is always $2$-dimensional, while we have proved that the Jacobian has always a $3$-dimensional kernel on $\mathcal{S}$.
\end{proof}

\begin{corollary}
	Magnetic Particle Imaging with a Field-Free Plane topology is not possible.
\end{corollary}

\begin{proof}
	Magnetic fields must obey Maxwell's equations and in the particular case of fields without magnetic field sources and without magnetization - as in the fields produced inside an MPI scanner - fulfill both magneto-static and quasi-static approximations~\cite{bringout_phdthesis}. This translates mathematically into the fact that the fields satisfy \eqref{eq:zero:div} and \eqref{eq:zero:rot}. A Field-Free Plane is a $2$-dimensional regular surface defined as the level set $\lbrace \bm{x}\in\Omega\colon\bm{H}(\bm{x})=0\rbrace$, which is impossible by Theorem \ref{thm:no:ffplanes}.
\end{proof}

In view of this result and of the fact that it is possible to move and rotate the FFL in any position in space~\cite{Knopp2010FFLformation,Soydan20233DtrajectoryFFL}, we will consider the FFL topology and propose model-based reconstruction formulae for 3D FFL MPI.

\section{Field-Free Line Magnetic Particle Imaging in 3D}\label{sec:3d:ffl}

In this section we introduce our proposed scanning setup for 3D FFL scans, the FFL MPI Core Operator and the results that allow us to formulate the sought reconstruction formulae.
We start by describing our proposed 3D FFL scanning modality which is based on the model employed to generate FFLs on a plane~\cite{Knopp2010FFLformation,Knopp_2011_Fourierslicethm}. Because we build our proposed scanning modality on the modeling of magnetic fields employed in the 2D FFL, we start by recalling its mathematical formulation.

Consider the frame of reference $(x,y,z)$ of the scanner. We will denote from now on with $\bm{e}_x$ the unitary vector in the direction of the $x$-axis and $\bm{e}_y$,$\bm{e}_z$ analogously. The static magnetic field with a FFL along the $x$-axis and gradient strength $G\in\mathbb{R}$ has the following expression~\cite{Knopp_2011_Fourierslicethm,Knopp2010FFLformation}:
\begin{displaymath}
	\bm{H}_S^0 (\bm{x} ) = G \begin{pmatrix}
		0 & 0 & 0 \\
		0 & -1 & 0 \\
		0 & 0 & 1
	\end{pmatrix}
	\bm{x} .
\end{displaymath}
Consider a rotation angle $\theta\in [0,2\pi )$, the magnetic field $\bm{H}_S^\theta$ obtained by rotating counterclockwise $\bm{H}_S^0$ in the $xy$-plane by $\theta$ can be written as:
\begin{equation}\label{eq:theta:field:x}
	\bm{H}_S^\theta (\bm{x} ) = G R_\theta \begin{pmatrix}
		0 & 0 & 0 \\
		0 & -1 & 0 \\
		0 & 0 & 1
	\end{pmatrix}
	R_\theta^T\bm{x} ,
\end{equation}
where $R_\theta$ is the rotation matrix
\begin{equation}\label{eq:rot:matrix}
	R_\theta = \begin{pmatrix}
		\cos\theta & -\sin\theta & 0\\
		\sin\theta & \cos\theta & 0 \\
		0 & 0 & 1
	\end{pmatrix}.
\end{equation}

Let now $(\bm{e}_\theta ,\bm{e}_\theta^\perp ,\bm{e}_z)$ be the orthonormal rotated basis consisting of the columns of $R_\theta$, i.e., $\bm{e}_\theta = (\cos\theta ,\sin\theta , 0)^T$, $\bm{e}_\theta^\perp = (-\sin\theta ,\cos\theta , 0)^T$ and $\bm{e}_z$. We can write the components of the static field $\bm{H}_S^\theta$ with respect to the rotated basis by simply performing the matrix multiplications in \eqref{eq:theta:field:x}:
\begin{equation}
	\begin{split}
		\bm{H}_S^\theta (\bm{x}) & = G (x\sin\theta - y\cos\theta )\bm{e}_\theta^\perp + Gz\bm{e}_z \\
		& = -G \langle \bm{x}\, ,\bm{e}_\theta^\perp\rangle \bm{e}_\theta^\perp + G\langle \bm{x}\, , \bm{e}_z\rangle\bm{e}_z ,
	\end{split}
\end{equation}
where $\langle\cdot\, , \cdot\rangle$ is the standard scalar product. In the 2D scenario, to move the FFL in the $xy$-plane while keeping it parallel to the $\bm{e}_\theta$ direction, an oscillating dynamic field $\bm{H}_D^\theta$ along the orthogonal direction $\bm{e}_\theta^\perp$ is added:
\begin{equation}\label{eq:2dffl:dynamic:field}
	\bm{H}_D^\theta (t) = A \Lambda (t)\bm{e}_\theta^\perp ,
\end{equation}
with $\Lambda (t)$ usually chosen to be sinusoidal~\cite{Knopp_2011_Fourierslicethm} of amplitude $A$, e.g.,  $\Lambda (t) = \cos (2\pi f_0 t)$ for some excitation frequency $f_0$. Hence, the magnetic field applied for a scan given an angle $\theta$ is the superposition of the static and dynamic fields:
\begin{equation}\label{eq:2dffl:field}
	\bm{H}^\theta (\bm{x},t) = \bm{H}_S^\theta (\bm{x}) + \bm{H}_D^\theta (t) .
\end{equation}

\subsection{Proposed Scanning Setup}\label{subsec:proposed:ffl3d}
We now describe a possible purely 3D FFL scanning modality. Our idea is to consider the 2D FFL setup in the $xy$-plane described in \eqref{eq:2dffl:field} with dynamic field applied as in \eqref{eq:2dffl:dynamic:field} and allow movement of the FFL also along the $z$-axis (this is possible in view of the recent results in~\cite{Soydan20233DtrajectoryFFL,vogel2014TWMPI,vogel2023impi}). Therefore, given an angle $\theta$ and a static magnetic field $\bm{H}_S^\theta (\bm{x})$ as in \eqref{eq:theta:field:x}, we consider the following dynamic field:
\begin{equation}\label{eq:dynamic:field}
	\bm{H}_D^\theta (t) = A_1\Lambda_1 (t)\bm{e}_\theta^\perp - A_2\Lambda_2 (t)\bm{e}_z ,
\end{equation}
with sinusoidal functions $\Lambda_1$,$\Lambda_2$ and amplitudes $A_1$,$A_2$. The magnetic field applied to the particle distribution in this 3D FFL scan is therefore
\begin{equation}\label{eq:magfield}
	\begin{split}
		\bm{H}^\theta (\bm{x},t) & = \bm{H}_S^\theta (\bm{x})+\bm{H}_D^\theta (t) \\
		& = \left ( A_1\Lambda_1 (t)-G\langle\bm{x},\bm{e}_\theta^\perp\rangle\right ) \bm{e}_\theta^\perp + \left ( -A_2\Lambda_2 (t)+Gz\right ) \bm{e}_z .
	\end{split}
\end{equation}
Using \eqref{eq:faraday} and \eqref{eq:signal:encoding} and assuming that the receive coil sensitivity pattern is constant,i.e., $P(\bm{x})=P\in\mathbb{R}^{3\times 3}$, the signal obtained during a scan applying the magnetic field $\bm{H}^\theta (\bm{x},t)$ is given by
\begin{equation}\label{eq:signal:theta}
	\bm{s}^\theta (t) = -\mu_0 m P \frac{d}{dt}\int_{\mathbb{R}^3}\rho (\bm{x})\mathcal{L}\biggl (\frac{\lVert\bm{H}^\theta (\bm{x},t)\rVert}{H_{\mathrm{sat}}}\biggr )\frac{\bm{H}^\theta (\bm{x},t)}{\lVert\bm{H}^\theta (\bm{x},t)\rVert}\, d\bm{x} .
\end{equation}

\subsection{Decomposition of the Signal and the MPI Core Operator in the FFL setting}

In this section we state and prove the main results of this work (Theorem \ref{thm:main}). In particular, we show that for the 3D FFL scan proposed with the magnetic fields as in \eqref{eq:magfield} the signal obtained is related to the distribution $\rho$ via the MPI Core Operator and the X-Ray projection operator. Before stating the theorem we recall some important mathematical tools needed for its formulation.

The first fundamental tool we recall is the matrix-valued kernel $K_h$ and the MPI Core Operator~\cite{marz2016model}, of which we give here a definition.
\begin{definition}\label{def:mpi:tools}
	For any $n\geq 1$, $h>0$ let $K_h\colon\mathbb{R}^n\to\mathbb{R}^{n\times n}$ the rotationally invariant matrix-valued kernel defined as 
	\begin{equation}\label{eq:mpi:kernel}
		K_h (\bm{y}) = 
		\nabla_{\bm{y}} \left ( \mathcal{L}\left ( \frac{\lVert \bm{y}\rVert}{h}\right )\frac{\bm{y}}{\lVert\bm{y}\rVert}\right ) =
		\frac{1}{h} \left [ f_1 \left ( \frac{\lVert \bm{y}\rVert}{h}\right ) \frac{\bm{y}\bm{y}^T}{\lVert\bm{y}\rVert^2} + f_2\left (\frac{\lVert\bm{y}\rVert}{h}\right )\left ( \mathbb{I}_n - \frac{\bm{y}\bm{y}^T}{\lVert\bm{y}\rVert^2}\right ) \right ] ,
	\end{equation}
	where $f_1 (z)=\mathcal{L}'(z)$ and $f_2 (z) = \mathcal{L} (z) / z$ for $z\in\mathbb{R}$, extended in zero by continuity.
	The MPI Core Operator is the bounded linear operator $A_h\colon L^1 (\mathbb{R}^n )\to C^{\infty}_b (\mathbb{R}^n ;\mathbb{R}^{n\times n})$ for $f\in L^1 (\mathbb{R}^n)$ and $C^{\infty}_b (\mathbb{R}^n ;\mathbb{R}^{n\times n})$ the set of smooth and bounded functions as
	\begin{equation}\label{eq:mpi:core:op}
		A_h [f] (\bm{r}) \coloneq \left ( K_h * f \right )(\bm{r}) = \int_{\mathbb{R}^n} f(\bm{x}) K_h (\bm{r}-\bm{x})\, d\bm{x} ,
	\end{equation}
	where the convolution is intended to be componentwise and $A_h$ is well-defined.
\end{definition}

\begin{rmk}
	The functions $f_1 (\lVert\cdot\rVert)$ and $f_2 (\lVert\cdot\rVert)$ are both bounded and in $L^p (\mathbb{R}^n )$ if $p>n$. It follows that $K_h\in L^{n+\varepsilon}(\mathbb{R}^n ;\mathbb{R}^{n\times n})$ for every $\varepsilon >0$ and it follows from Young's convolution inequality that $A_h\colon L^{p}(\mathbb{R}^n )\to L^q (\mathbb{R}^n ;\mathbb{R}^{n\times n})$ is a bounded linear operator for every $p,q\geq 1$ such that $1+\frac{1}{q}=\frac{1}{p} + \frac{1}{n+\varepsilon}$.
\end{rmk}

The MPI Core Operator defined in \eqref{eq:mpi:core:op} is the integral operator that is obtained by passing the time derivative under the integral sign in \eqref{eq:signal:theta} and performing suitable changes of variables (we refer the interested reader to~\cite{marz2016model} for the details).

We will show how this operator encapsulates the relationship between the signal, the data describing the scans and the distribution of particles. In what follows, the MPI Core Operator will not act directly on the target distribution $\rho$ but on its X-Ray projection~\cite{helgason2011,NattererCTbook}, which we briefly recall. 
If $\Omega\subset\mathbb{R}^n$ with $n\geq 3$ is an non-empty, open and bounded set, the X-Ray transform is the operator $X\colon L^1 (\Omega )\times T\mathbb{S}^{n-1}\to \mathbb{R}$ defined as follows\footnote{Here  $T\mathbb{S}^{n-1}$ denotes the tangent bundle of the unit $n$-sphere $\mathbb{S}^{n-1}$.}:
\begin{equation*}
	(\rho ,\bm{\omega},\bm{x})\mapsto X\rho (\bm{\omega} ,\bm{x})=\int_{-\infty}^{+\infty}\rho (\bm{x}+s\bm{\omega})\, ds\quad\text{for}\quad\bm{\omega}\in\mathbb{S}^{n-1},\,\bm{x}\in\bm{\Pi}_{\bm{\omega}}
\end{equation*}
where $\bm{\Pi}_{\bm{\omega}}=\langle \bm{\omega}\rangle^\perp$ is the orthogonal space to the line spanned by $\bm{\omega}$, which we denote with the symbol $\langle \bm{\omega}\rangle$. The operator X defines the projection operator $X_{\bm{\omega}}[\rho ]$ onto the orthogonal space $\bm{\Pi}_{\bm{\omega}}$ in the following way:
\begin{equation*}
	X_{\bm{\omega}}[\rho ]\colon \bm{\Pi}_{\bm{\omega}}\to \mathbb{R}\quad ,\quad \bm{x}\mapsto X_{\bm{\omega}}[\rho ](\bm{x})\coloneq X\rho (\bm{\omega},\bm{x}) .
\end{equation*}
Finally, given an angle $\theta$, we define the linear transformation $E_\theta\coloneq\diag (1,1,-1)\cdot R_\theta$ for $R_\theta$ in \eqref{eq:rot:matrix}.

With the tools introduced we can state the results that links the signal to the target distribution and allows us to deduce reconstruction formulae for 3D FFL scans for functions in fractional Sobolev spaces, whose definition we now recall from \cite{bony2001cours}. We denote with $\mathcal{S}'(\mathbb{R}^n )$ the space of tempered distributions and given $\alpha > 0 $, we consider $$H^{\alpha } (\mathbb{R}^n ) = \left\lbrace f\in\mathcal{S}'(\mathbb{R}^n) \colon \int_{\mathbb{R}^n}(1+\lvert\xi\rvert^2 )^\alpha \lvert \Fourier f(\xi )\rvert^2\, d\xi ) <\infty\right\rbrace .$$ For $\Omega$ open, $H^{\alpha }(\Omega )$ is the set of functions $f$ for which there exists a $g\in H^{\alpha }(\mathbb{R}^n )$ extending $f$.
\begin{theorem}\label{thm:main}
	Let $\Omega\subset\mathbb{R}^3$ be a non-empty, open set with Lipschitz and bounded boundary, $\rho\in H^{s+1} (\Omega )$ with $s>0$ and $\rho = 0$ on $\mathbb{R}^{3}\setminus \overline{\Omega}$ a particle distribution, $P\in\mathbb{R}^{3\times 3}$ a homogeneous receive coil sensitivity, $\theta\in [0\, ,\pi )$ a scan angle and $\bm{H}^\theta (\bm{x},t)$ a dynamic magnetic field producing a FFL line as in \eqref{eq:magfield}, the induced voltage $\bm{s}^\theta (t)$ collected by the receiving coil is related to the distribution $\rho$ via the MPI Core Operator and the X-Ray transform by the equation
	\begin{equation}\label{eq:main}
		\bm{s}^\theta (t) = -\mu_0 m PE_\theta A_{h}\bigl [X_{\bm{e}_\theta} [\rho]\bigr ]\bigl (\bm{r}(t)\bigr )\bm{v}(t) ,
	\end{equation}
	where $\bm{r}(t)$ is the trajectory of the point where the FFL  and the orthogonal plane $\bm{\Pi}_\theta=\langle (\cos\theta \, ,\sin\theta \, , 0)^T\rangle^\perp$ intersect and $\bm{v}(t)=\frac{d}{dt}\bm{r}(t)$.
\end{theorem}

\begin{proof}
	Consider the change of variables $\bm{\zeta}=(\eta ,\xi ,z)=R_\theta^T\bm{x}$, i.e., 
	\begin{equation}\label{eq:change:var}
		\begin{cases}
			\eta  = x\cos\theta +y\sin\theta \\
			\xi  = -x\sin\theta +y\cos\theta \\
			z  = z 
		\end{cases},
	\end{equation}
	then $\bm{x}=R_\theta\bm{\eta}$ and the magnetic field can be rewritten as
	\begin{align}\label{eq:b:eta}
		\bm{H}^\theta (\bm{x},t) & = \bigl ( A_1\Lambda_1 (t)-G\langle R_\theta\bm{\eta},\bm{e}_\theta^\perp\rangle\bigr ) \bm{e}_\theta^\perp + \bigl (-A_2\Lambda_2 (t)+Gz\bigr )\bm{e}_z \notag\\
		& = \bigl ( A_1\Lambda_1 (t)-G\langle \bm{\eta},R_\theta^T\bm{e}_\theta^\perp\rangle\bigr ) \bm{e}_\theta^\perp + \bigl (-A_2\Lambda_2 (t)+Gz\bigr )\bm{e}_z \notag \\
		& = \bigl ( A_1\Lambda_1 (t)-G\xi\bigr ) \bm{e}_\theta^\perp + \bigl (-A_2\Lambda_2 (t)+Gz\bigr )\bm{e}_z \notag \\
		& = E_\theta \hat{\bm{H}}^\theta (\xi ,z,t) ,
	\end{align}
	where 
	\begin{equation}\label{eq:b:hat}
		E_\theta = \begin{pmatrix}
			\cos\theta & -\sin\theta & 0\\
			\sin\theta & \cos\theta & 0 \\
			0 & 0 & -1
		\end{pmatrix}\quad\text{and}\quad\hat{\bm{H}}^\theta (\xi ,z,t) = \begin{pmatrix}
			0 \\ A_1\Lambda_1 (t)-G\xi \\ A_2 \Lambda_2 (t)-Gz
		\end{pmatrix}.
	\end{equation}
	We now use change of variable - with $\det R_\theta = 1$ - and substitute the expression in \eqref{eq:b:eta} into the signal in \eqref{eq:signal:theta}. The magnetic field $\hat{\bm{H}}^\theta (\xi ,z,t)$ is independent of the variable $\eta$. Thus, using Fubini-Tonelli's theorem we can integrate separately $\rho$ along the FFL, i.e., along the direction $\bm{e}_\eta = \bm{e}_\theta$ to obtain
	\begin{align}\label{eq:signal:rot:var}
		\bm{s}^\theta (t) & = -\mu_0 m P E_\theta \frac{d}{dt}\int_{\mathbb{R}^3}\rho (\eta,\xi,z)\mathcal{L}\biggl (\frac{\lVert\hat{\bm{H}}^\theta (\xi,z,t)\rVert}{H_{\mathrm{sat}}}\biggr )\frac{\hat{\bm{H}}^\theta (\xi,z,t)}{\lVert\hat{\bm{H}}^\theta (\xi,z,t)\rVert}\, d\bm{\zeta} \notag\\
		& =  -\mu_0 m P E_\theta \frac{d}{dt}\int_{\mathbb{R}^2}\underbrace{\biggl (\int_\mathbb{R}\rho (\eta ,\xi ,z)\, d\eta\biggr )}_{\eqqcolon \hat{\rho}_\theta (\xi ,z)}\mathcal{L}\biggl (\frac{\lVert\hat{\bm{H}}^\theta (\xi,z,t)\rVert}{H_{\mathrm{sat}}}\biggr )\frac{\hat{\bm{H}}^\theta (\xi,z,t)}{\lVert\hat{\bm{H}}^\theta (\xi,z,t)\rVert}\, d\xi dz \notag \\
		& = -\mu_0 m P E_\theta \frac{d}{dt}\int_{\mathbb{R}^2} \hat{\rho}_\theta (\xi ,z)\mathcal{L}\biggl (\frac{\lVert\hat{\bm{H}}^\theta (\xi,z,t)\rVert}{H_{\mathrm{sat}}}\biggr )\frac{\hat{\bm{H}}^\theta (\xi,z,t)}{\lVert\hat{\bm{H}}^\theta (\xi,z,t)\rVert}\, d\xi dz ,
	\end{align}
	where $\hat{\rho}_{\theta} (\xi ,z) = X_{\bm{e}_\theta} [\rho ](\xi ,z)$ is the X-Ray transform of $\rho$ in the direction $\bm{\bm{e}_\theta}$. The X-Ray projection in our setting is well-defined using the trace theorem; in particular, let $\bm{L}$ be a line and $\bm{\Pi}_{\bm{L}}\supset\bm{L}$ a plane, one has the following applications of the (surjective) trace operators $H^{s+1}(\Omega)\xrightarrow{\trace_1} H^{s+\frac{1}{2}}(\Omega\cap \bm{\Pi}_{\bm{L}} )\xrightarrow{\trace_2}H^{s}(\Omega\cap \bm{L}) \xhookrightarrow{\mathrm{\iota}}L^{1}(\Omega\cap \bm{L}) $ where the last inclusion $\iota$ is true because $\Omega$ is bounded. Consequently, one can define the integral along any line by considering the composition of trace operators.
	Now, if we multiply and divide the magnetic field in \eqref{eq:signal:rot:var} by the gradient strength $G$ and rename the variables suitably, i.e., we consider
	\begin{displaymath}
		\frac{\hat{\bm{H}}^\theta (\xi ,z ,t)}{G} = \biggl (		0 \, , \frac{A_1\Lambda_1 (t)}{G}-\xi \, , \frac{A_2\Lambda_2 (t)}{G}-z \biggr )^T
	\end{displaymath}
	and label the variables with $\bm{y}=(0,\xi ,z)$ on the plane $\bm{\Pi}_\theta =\langle \bm{e}_\theta\rangle^\perp$, with
	\begin{equation}\label{eq:lissajous:plane}
		\bm{r}(t) = \biggl (		0 \, , \frac{A_1\Lambda_1 (t)}{G} \, , \frac{A_2\Lambda_2 (t)}{G} \biggr )^T ,
	\end{equation}
	the trajectory of the intersection point between the FFL at time $t$ and the plane $\bm{\Pi}_\theta$, we can rewrite the signal as
	\begin{align}
		\bm{s}^\theta (t) & = -\mu_0 mP E_\theta \frac{d}{dt}\int_{\mathbb{R}^2}\hat{\rho}_\theta (\bm{y}) \mathcal{L}\biggl ( \frac{\lVert \bm{r}(t)-\bm{y}\rVert}{h}\biggr ) \frac{ \bm{r}(t)-\bm{y}}{\lVert \bm{r}(t)-\bm{y}\rVert}\, d\bm{y} \notag\\
		& = -\mu_0 mP E_\theta \biggl [\int_{\mathbb{R}^2}\hat{\rho}_\theta (\bm{y}) \nabla_{\bm{r}}\biggl (\mathcal{L}\biggl ( \frac{\lVert \bm{r}(t)-\bm{y}\rVert}{h}\biggr ) \frac{ \bm{r}(t)-\bm{y}}{\lVert \bm{r}(t)-\bm{y}\rVert}\biggr )\, d\bm{y} \biggr ]\bm{v}(t) ,\label{eq:signal:conv:2}
	\end{align}
	where $h=\frac{H_{\mathrm{sat}}}{G}$, $\bm{v}(t)=\frac{d\bm{r}(t)}{dt}$ and we have differentiated under the integral to get \eqref{eq:signal:conv:2}.
	
	Using Definition \ref{def:mpi:tools} together with \eqref{eq:signal:conv:2}, the relationship between the data $\bm{s}^\theta (t)$ and the distribution $\rho$ is mediated by the MPI Core Operator and the X-Ray transform by the equation
	\begin{displaymath}
		\bm{s}^\theta (t) = -\mu_0 m PE_\theta A_{h}\bigl [X_{\bm{e}_\theta} [\rho]\bigr ]\bigl (\bm{r}(t)\bigr )\bm{v}(t) ,
	\end{displaymath}
	which concludes the proof of the theorem.
\end{proof}

\begin{rmk}\label{rmk:reconstruction}
	We assume that blood vessel structures are reasonably approximated by a compact set $K=\overline{\Omega}\subset\mathbb{R}^3$ with $\partial K = \partial\Omega $ of Lipschitz class. Consider then a function $g\in H^{s+1}(\Omega )$ and let $\rho = g\chi_K$ (extending $\rho$ onto $\partial\Omega$ with the trace operator); due to the possible discontinuity along $\partial\Omega$, $\rho=g\chi_K \not\in H^{s+1}(\mathbb{R}^{3})$ but Theorem \ref{thm:main} holds true. In Theorem \ref{thm:main} we show that the X-Ray projections of $\rho$ are still defined and in $L^{1}(\bm{\Pi}_{\theta})$, which means that $A_h\left  [X_{\bm{e}_\theta} [\rho]\right ]$ in Definition \ref{def:mpi:tools} is well defined and is a $C^{\infty}_b$ function. If then, $\bm{r}\colon [0,T]\to \mathbb{R}^3$ is a trajectory of class $C^k$, $k\geq 1$, then the theoretical signal $\bm{s}^{\theta} (t)$ in \eqref{eq:tr:main} is a $C^{k-1}$ function on the compact $[0,T]$.
\end{rmk}

We now clarify the relationship between the signal $\bm{s}^\theta (t)$ and the target distribution $\rho$ made explicit in Theorem \ref{thm:main}. In particular, given an angle $\theta$, the intersection of the FFL and the plane $\bm{\Pi}_\theta$ orthogonal to the direction $\bm{e}_\theta$, i.e., $\bm{\Pi}_\theta =\langle \bm{e}_\theta\rangle^\perp$, we have in our setup that the intersection between the FFL and the plane $\bm{\Pi}_\theta$ draws in time a path $t\mapsto\bm{r}(t)\in\bm{\Pi}_\theta$ (see Fig.\ref{fig:scan}) with velocity $\bm{v}(t)$; the curve $\bm{r}(t)$ is the point of evaluation of the MPI Core Operator $A_{h}$ in \eqref{eq:main}, which is applied to the X-Ray projection $X_{\bm{e}_\theta} [\rho ]$ in the direction $\bm{e}_\theta$ of the FFL. Remarkably, the relationship in \eqref{eq:main} bears similarities with the decomposition of the signal coming from a FFP 2D scans~\cite{marz2016model}. Indeed, in 3D FFP it is possible~\cite{marz2016model} to rewrite the magnetic field $\bm{H}(\bm{x},t)$ only in terms of the gradient matrix $\bm{G}=G\cdot \mathrm{diag}(-1,-1,2)$ and the trajectory of the FFP $\bm{r}(t)$ as $\bm{H}(\bm{x},t)=\bm{G}(\bm{x}-\bm{r}(t))$. From these computations it follows that the induced voltage $\bm{s}(t)$ in \eqref{eq:signal:encoding} produced by a FFP scan along the trajectory $\bm{r}(t)$ is
\begin{equation}\label{eq:rel:ffp}
	\bm{s}(t)=\mu_0 m \biggl [P \int_{\mathbb{R}^3}\rho (\bm{x})\nabla_{\bm{r}}\biggl (\mathcal{L}\biggl (\frac{\lvert \bm{G}(\bm{r}(t)-\bm{x})\rvert}{H_{\mathrm{sat}}}\biggr )\frac{ \bm{G}(\bm{r}(t)-\bm{x})}{\lvert \bm{G}(\bm{r}(t)-\bm{x})\rvert}\biggr )\, d\bm{x}\biggr ]\bm{v}(t) ,
\end{equation}
where $P=P(\bm{x})$ is the sensitivity pattern of the receive coils that can be assumed homogeneous and $\bm{v}(t)$=$\frac{d}{dt}\bm{r}(t)$.

Given Definition \ref{def:mpi:tools} of the MPI Core Operator and performing change of variables (see~\cite{marz2016model} for details) that get rid of the terms $\bm{G}$,$P$,$m$ and $\mu_0$ in \eqref{eq:rel:ffp}, the data $\bm{s}(t)$ produced during a 3D FFP scan can be rewritten as
\begin{equation}\label{eq:ffp:relation}
	\bm{s}(t)= A_{h}[\rho ]\bigl (\bm{r}(t)\bigr )\bm{v}(t) .
\end{equation}
\eqref{eq:ffp:relation} holds also true for 1D and 2D FFP scans~\cite{marz2016model}, with $\bm{r}(t)$ and $\bm{v}(t)$ moving on a line or inside a plane, respectively.
\eqref{eq:ffp:relation} is incredibly similar to the relation found in \eqref{eq:tr:main}, which can be thought of as a 2D FFP scan of the X-Ray projection, with the terms $P$ and $E_\theta$ that describe the change of the signal induced in the receiving coil depending on the angle $\theta$ considered. This similarity is graphically illustrated in Fig.\ref{fig:scan}.

\def\imratio{0.6}
\begin{figure}[t]
	\centering
	\begin{subfigure}[t]{\imratio\linewidth}
		\includegraphics[width=\linewidth]{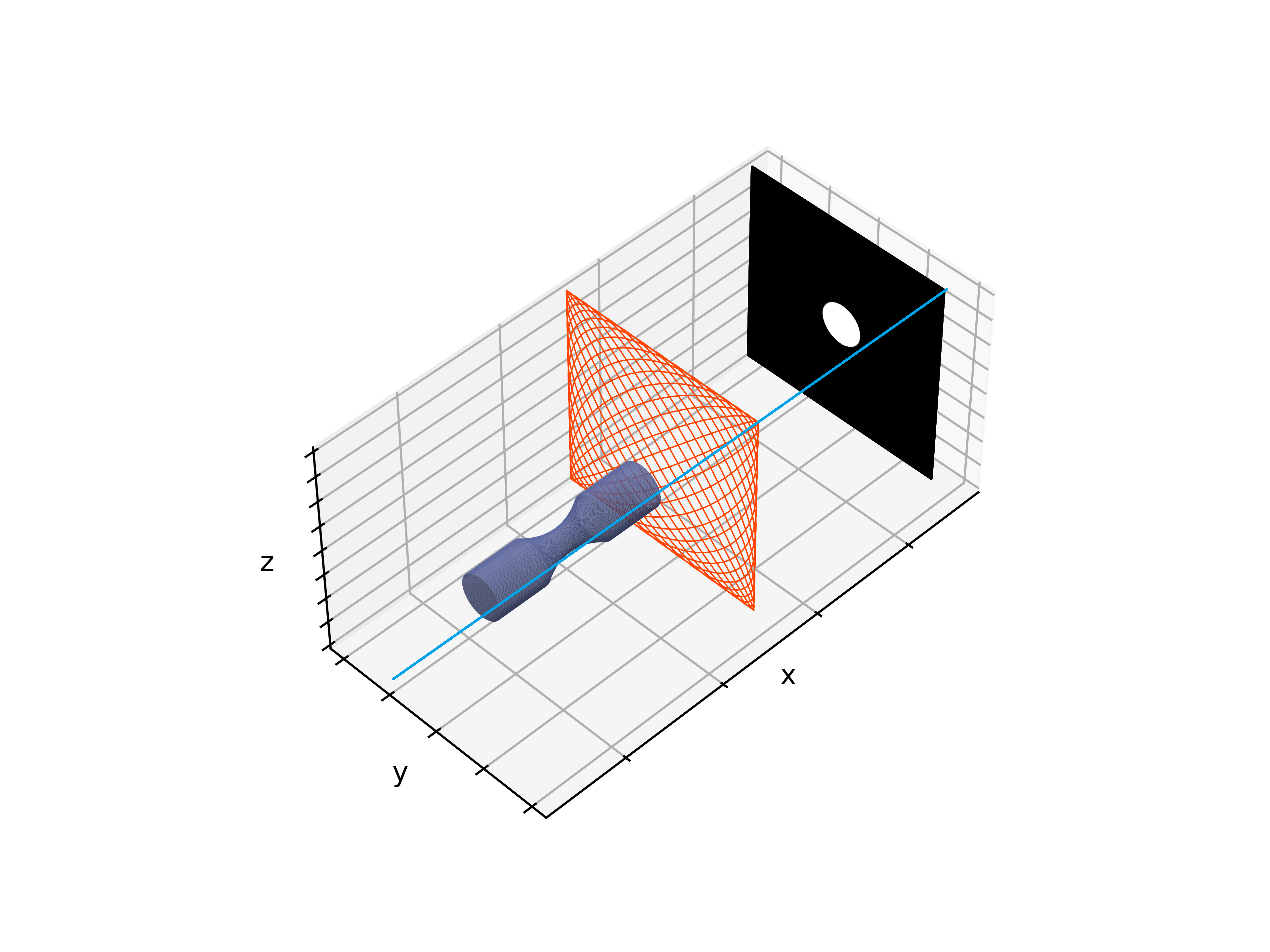}
	\end{subfigure}
	\caption{A graphic representation of the components of the relation in \eqref{eq:main} for $\theta = 0$. The vessel-shaped particle distribution $\rho$ (in blue) is scanned moving the FFL (light blue) and the intersection point $\bm{r}(t)$ between the FFL and the $yz$-axis forms a Lissajous trajectory (in orange). The picture in black and white orthogonal to the $x$-axis is the X-Ray projection $X_{\bm{e}_\theta} [\rho ]$ for $\theta = 0$. }
	\label{fig:scan}
\end{figure}

\subsection{Reconstruction formulae}

In \eqref{eq:main} we have learned that the target distribution $\rho$ is linked to the signal $\bm{s}^\theta$ via the MPI Core Operator and consequently, the reconstruction of $\rho$ from the data must pass through the inversion of the MPI Core Operator. The inversion of the MPI Core Operator has been explored in~\cite{marz2016model}, to which we refer the interested reader for a detailed discussion of the formulas and results presented in this section and in particular for a discussion on the ill-posedness of the inverse problem associated to the inversion. Before discussing further, we transform the data $\bm{s}^\theta (t)$ for every $\theta$ in a way such that the relationship in \eqref{eq:main} can be rewritten as the result of a direct application of the MPI Core Operator and define
\begin{equation}\label{eq:tr:main}
	\tilde{\bm{s}}^\theta (t) = A_{h}\bigl [X_{\bm{e}_\theta} [\rho]\bigr ]\bigl (\bm{r}(t)\bigr )\bm{v}(t) \quad\text{with}\quad \tilde{\bm{s}}^\theta (t) = -\frac{1}{m\mu_0 }E_\theta^{-1} P^{-1}\bm{s}^\theta .
\end{equation}
The relation in \eqref{eq:tr:main} is precisely of the form studied in~\cite{marz2016model} and it has been proven that the trace of the MPI Core Operator contains all the information needed for the reconstruction, i.e., if $\alpha_{h}[\rho ](r) = \trace A_h [\rho ]$ for some distribution $\rho$, then the trace is the convolution $\alpha_{h}[\rho ](r) = (\rho * \kappa_h )(r)$ between the distribution $\rho$ and the scalar valued kernel $\kappa_h$ having the form
\begin{equation}\label{eq:conv:kernel}
	\kappa_h (y) = \frac{1}{h}f\left ( \frac{\lvert y\rvert}{h}\right ) , \quad\text{with}\quad f(z)=\mathcal{L}' (z) + \mathcal{L} (z)\frac{n-1}{z} 
\end{equation}
in the general n-dimensional case.
\begin{rmk}\label{rmk:injective}
	We remark that the convolution operator $\alpha_h [f] \coloneq \kappa_h * f$ is an injective operator from $L^1 (\mathbb{R}^2)$ functions with compact support into $C^{\infty}_b(\mathbb{R}^2)$. Indeed, $\kappa_h\in \mathcal{S}' (\mathbb{R}^2)$ and $f\in L^1$ and compactly supported (hence in $ \mathcal{E}' (\mathbb{R}^2)$) and consequently, the Fourier convolution theorem implies that $\alpha_h$ is injective if and only if $\Fourier \kappa_h$ has no zeros. This is indeed the case using recent analysis \cite{maass2019fourier} that allow to show that $\Fourier\kappa_h (\bm{\xi}) = -2\pi \mathcal{H}'_0 [f_2 ](\lvert\bm{\xi}\rvert ) \lvert\bm{\xi}\rvert$ and the derivative of the Hankel transform $\mathcal{H}_0$ is strictly negative with exponential singularity in 0.
\end{rmk}
\noindent
With these considerations we can state the first reconstruction formula:
\begin{theorem}[First Reconstruction Formula]\label{thm:first:formula}
	Let $\rho$ be as in Theorem \ref{thm:main}, $\theta\in [0,\pi)$ and $F(r)=A_{h} \bigl [ X_{\bm{e}_\theta} [\rho ]\bigr ](r)$ be an ideal sampling of the MPI Core Operator each point $r$, then 
	\begin{equation}\label{eq:first:formula}
		X_{\bm{e}_\theta} [\rho ] = \kappa_h^{-1}\circ \trace \left (F\right ) ,
	\end{equation}
	where we denote with the symbol $\kappa_h^{-1}\circ$ the deconvolution w.r.t. the kernel $\kappa_h$.
\end{theorem}
\begin{proof}
	Let $\theta\in [0,\pi)$ be a scanning angle, the signal $\bm{s}^{\theta}(t)$ produced by a scan can be described according to Theorem \ref{thm:main} as in \eqref{eq:signal:conv:2}. Applying the change of variables in \eqref{eq:tr:main}, we obtain that the relationship between $\tilde{\bm{s}}^\theta (t)$ and the X-Ray projection $X_{\bm{e}_\theta} [\rho]$ can be written as
	\begin{equation}\label{eq:ffl:2d:data}
		\tilde{\bm{s}}^\theta (t) = \biggl [\int_{\mathbb{R}^2}X_{\bm{e}_\theta} [\rho] (\bm{y}) \nabla_{\bm{r}}\biggl (\mathcal{L}\biggl ( \frac{\lVert \bm{r}(t)-\bm{y}\rVert}{h}\biggr ) \frac{ \bm{r}(t)-\bm{y}}{\lVert \bm{r}(t)-\bm{y}\rVert}\biggr )\, d\bm{y} \biggr ]\bm{v}(t) ,
	\end{equation}
	with $\bm{r}(t)$ as in \eqref{eq:lissajous:plane}, and $\bm{v}(t)=\frac{d}{dt}\bm{r}(t)$. In particular, $\bm{r}(t)$ as in \eqref{eq:lissajous:plane} is a 2D trajectory in the $yz$-plane (the first coordinate is always $0$), the X-Ray projection $X_{\bm{e}_\theta} [\rho ]$ is a distribution in 2 variables and the integral in \eqref{eq:ffl:2d:data} is the 2D FFP MPI Core Operator as in \eqref{eq:rel:ffp}. This means that upon the change of variables, the modeling of the FFL 3D scan data is formally equivalent to the data obtained by performing a 2D FFP scan of the X-Ray projection as 2D distribution in the $yz$-plane. By Remark \ref{rmk:injective}, $\kappa_h$ is injective and therefore is a bijection onto its range, but $F\in\mathrm{range}(\kappa_h )$ by assumption of ideal sampling.
\end{proof}
\begin{rmk}
	The first reconstruction formula in Theorem \ref{thm:first:formula} involves a deconvolution which has been proven~\cite{marz2016model} to be ill-posed and constitutes the first source of ill-posedness of the reconstruction algorithm in Section \ref{sec:rec:algorithm}.
\end{rmk}

Applying the reconstruction formula in Theorem \ref{thm:first:formula}, one obtains a distribution $X_{\bm{e}_\theta}$ in two variables on the plane orthogonal to the direction $\bm{e}_\theta$ for every angle $\theta$. In particular, we obtain the sheaf of all planes through the $z$-axis and on each one of these planes we have the X-Ray projection of the distribution $\rho$. The idea behind the second reconstruction formula comes from the observation that the X-Ray and Radon transforms coincide in 2D~\cite{NattererCTbook}: indeed, if we restrict ourselves to a plane orthogonal to the $z$-axis, for example on the plane $\bm{\Pi}_a = \lbrace (x,y,z)\in\mathbb{R}^3 \colon z=a\rbrace$ for some $a\in\mathbb{R}$, then the 2D Radon transform of angle $\theta$ on $\bm{\Pi}_a$ coincides with the X-Ray transform of angle $\theta$ on $\bm{\Pi}_a$ and we can consider the Filtered Back Projection (FBP)~\cite{NattererCTbook}.
\begin{theorem}[Second Reconstruction Formula]\label{thm:second:formula}
	Given a 3D distribution $\rho$ as in Theorem \ref{thm:main} and the set of all of its X-Ray projections $X_{\bm{e}_\theta} [\rho ]$ for all $\theta\in [0,\pi )$ around the $z$-axis, then for every restriction of $\rho$ on the plane $\bm{\Pi}_a$ with $a\in\mathbb{R}$
	
	\begin{equation}\label{eq:second:formula}
		\rho (\bm{x} )= \frac{1}{4\pi} \mathcal{R}^* \left ( \mathcal{I}^{-1} X_{\bm{e}_\theta} [\rho ] \right ) (\bm{x}) \quad\text{for every }\bm{x}\in\bm{\Pi}_a 
	\end{equation}
	where $\mathcal{I}^{-1}$ is the Riesz potential \cite{NattererCTbook} and $\mathcal{R}^*$ is the adjoint operator (back projection) of the 2D Radon transform on $\bm{\Pi}_a$.
\end{theorem}
\begin{proof}
	Given $\bm{\Pi}_a = \lbrace (x,y,z)\in\mathbb{R}^3 \colon z=a\rbrace$ for some $a\in\mathbb{R}$, the restriction of the X-Ray projections $X_{\bm{e}_\theta} [\rho ]$ onto $\bm{\Pi}_a$ is the 2D Radon transform in the direction $\bm{e}_\theta = (\cos\theta\, ,\sin\theta )^T$ of the 2D distribution $\rho (\bm{x} )\bigl\lvert_{\bm{\Pi}_a}$, where both the restriction and the X-Ray transforms are well defined and compactly supported in view of Remark \ref{rmk:reconstruction}. The X-Ray projections lay on the planes $\bm{\Pi}_\theta$ forming a sheaf around the $z$-axis, meaning that the set of restrictions $X_{\bm{e}_\theta} [\rho ]\bigl\lvert_{\bm{\Pi}_a}$ for $\theta\in[0,\pi )$ are in $L^2$, they form a complete set of 2D Radon transforms for all directions and the Radon inversion formula in~\cite[Theorem 2.1]{NattererCTbook} holds true for the Riesz potential $\mathcal{I}^\alpha$ with $\alpha =0$.
\end{proof}

A practical implementation in applied scenarios of \eqref{eq:second:formula}, involves a windows function to filter the spectrum of $\mathcal{I}^{-1}X_{\bm{e}_\theta}$, i.e., one employs the regularized equation
\begin{equation}
	\rho (\bm{x} )\biggl\lvert_{\bm{\Pi}_a} = \int_{0}^{\pi}\mathcal{F}^{-1}\left ( \mathcal{F}\biggl (X_{\bm{e}_\theta} [\rho ]\biggl\lvert_{\mathbf{L}_a}\biggr ) (\omega )\cdot \lvert\omega\rvert\cdot w (\lvert\omega\rvert )\right )\, d\theta .
\end{equation}

\section{Reconstruction algorithm}\label{sec:rec:algorithm}

In this section we describe the discrete data, and we formulate the reconstruction algorithm, which can be subdivided into three mains steps: (i) reconstruction of the MPI Core Operators from the scan data for each angle; (ii) regularized deconvolution of the traces of the MPI Core Operators reconstructed in the first step; (iii) reconstruction of the 3D distribution using the reconstruction formula in Theorem \ref{thm:second:formula}.

First, we describe the discretization employed in the experiments: we consider a 3D particle distribution $\rho\colon\mathbb{R}^3\to\mathbb{R}^+$ such that $\mathrm{supp}(\rho )$ is compactly contained in the scanning domain $\Omega = [x_1 ,x_2 ]\times [y_1 ,y_2 ]\times [z_1 ,z_2 ]$. This domain is discretized considering an $\tilde{N}_x \times \tilde{N}_y \times \tilde{N}_z$ grid and identifying the $(i,j,k)$-th cell of the grid with its center $(x_i ,y_j ,z_k )$ defined as
\begin{equation}\label{eq:disc}
	\begin{split}
		x_i & = x_1 + \left ( 0.5 + i\right )\Delta x  ,\quad \text{where}\quad \Delta x = (x_2 -x_1 )/N_x\text{  and  }i = 0,\cdots ,\tilde{N}_x - 1 \\
		y_j & = y_1 + \left ( 0.5 + j\right )\Delta y  ,\quad \text{where}\quad \Delta y = (y_2 -y_1 )/N_y\text{  and  }j = 0,\cdots ,\tilde{N}_y - 1 \\
		z_k & = z_1 + \left ( 0.5 + k\right )\Delta z  ,\quad \text{where}\quad \Delta z = (z_2 -z_1 )/N_z\text{  and  }k = 0,\cdots ,\tilde{N}_z - 1 .\\
	\end{split}
\end{equation}

For the acquisition of the scans, we consider the angles $\theta$ in some set of angles $\Theta$, usually $\theta_l = l\cdot \frac{\pi}{\lvert\Theta\rvert}$ where $l=0,\dots ,\lvert\Theta\rvert -1$ and $\lvert\Theta\rvert\in\mathbb{N}$ is the total number of angles considered. For each angle $\theta_l$ we consider an acquisition time $T_l >0$ of the scan and we discretize the time considering equidistant points $t_m = m\cdot \frac{T_m}{L_l}$ for $m=0,\dots ,L_l-1$ and $L_l$ is the number of time points of the scan performed with angle $\theta_l$. With the discretization of time we obtain a discretization of the scanning trajectories as in \eqref{eq:lissajous:plane},
$\bm{r}_m = \bm{r}(t_m )$, the velocities $\bm{v}_m = \bm{v}(t_m )$ and of the signal $\bm{s}_m^{\theta_l}=\bm{s}^{\theta_l}(t_m )$ for every $m=0,\dots ,L_l-1$ and $l=0,\dots ,\lvert\Theta\rvert -1$.
The task of the reconstruction algorithm is to reconstruct a discrete version of the distribution from the data on a grid $N_x \times N_y \times N_z$ with $N_i < \tilde{N}_i$, i.e., to obtain a discrete distribution $\rho_{\mathrm{rec}}\in\mathbb{R}^{N_x\times N_y\times N_z}$ with an algorithm in three steps which we now describe in detail.

\paragraph{\textbf{First Stage}: Obtain the Trace of the MPI Core Operator from the Data} In this first stage we aim at reconstructing the MPI Core Operators from the raw data, which for each $\theta_l\in\Theta$ consists of the set of triples $\lbrace (\bm{s}^{\theta_l}_m , \bm{r}_m , \bm{v}_m )\colon m=0,\dots ,L_l-1\rbrace$. We transform the data points according to \eqref{eq:tr:main}, such that they are all in the same frame of reference. In particular, the problem for each angle is of the form $\bm{s}=A(\bm{r} )\bm{v}$, where $A(\bm{r} )\in\mathbb{R}^{2\times 2}$ is a shorthand notation for the discrete matrix-valued field $A_{j,k}\in\mathbb{R}^{2\times 2\times N_y\times N_z}$ evaluated at the point $\bm{r}$ and which discretizes the MPI Core Operator $A_{h}$; here $N_y\times N_z$ is the discretization of the transformed region of interest $[y_1 , y_2 ]\times [z_1 , z_2 ]$, which upon the transformation in \eqref{eq:tr:main} is the same for each angle. To solve this problem we choose to employ the variational approach introduced in \cite{gapyak2022mdpi}. In particular, we retrieve the MPI Core Operator $A(\bm{r})$ on a $N_y\times N_x$ grid by solving the minimization problem
\begin{equation}\label{eq:fisrt:stage}
	A^{\theta_l} = \arg\min_{\hat{A}}\biggl\lbrace \frac{1}{L_l}\sum_{m=1}^{L_l}\lVert \bm{s}_m^{\theta_l} - I[\hat{A}](\bm{r}_m )\bm{v}_m \rVert_2^2 + \mu \lVert D\hat{A}\rVert_2^2\biggr\rbrace ,
\end{equation}
where $D$ is the matrix represents the discretization of the gradient by first forward differences and $I$ is the interpolation operator for a chosen interpolation scheme. We have shown in \cite{gapyak2022mdpi} that the Euler-Lagrange equations of the functional in \eqref{eq:fisrt:stage} are of the form $GA = b$, for $A$ and $b$ vectors and $G$ a symmetric positive definite matrix (and solvable with the CG method). Once all $\lvert\Theta\rvert$ problems of the form in \eqref{eq:fisrt:stage} have been solved, the discrete field $u_{j,k}^{l}=\trace A_{j,k}^{\theta_l}$ of the traces of the reconstructed MPI Core Operator serves as input for the second stage.

\paragraph{\textbf{Second Stage}: Reconstruction of the X-Ray Projections from the Traces} The input of this stage are the trace fields $u^l$ for each angle $\theta_l\in\Theta$. From these trace fields it is possible to retrieve the X-Ray transforms of the distribution using the First Reconstruction Formula in Theorem \ref{thm:first:formula}. In particular, for each angle $\theta_l$ we solve the deconvolution problem in \eqref{eq:first:formula} of the trace field $u^l$ and the kernel $\kappa_h$ in \eqref{eq:conv:kernel} to obtain $\chi_l \coloneq X_{\theta_l}[\rho ]$. Because of course the data $\bm{s}^{\theta_l}$ and hence $u^l$ is corrupted by noise, the deconvolution problems are severely ill-posed and regularization is needed. Following~\cite{marz2016model} we employ the Tikhonov regularization technique with a smoothing regularizer (smoothness of the solution has been proven in~\cite{marz2016model}) and we formulate the (continuous) minimization problem
\begin{equation}\label{eq:second:cont}
	\chi_l = \arg\min_{\hat{\chi}} \lambda_l \left\lVert \nabla_r \hat{\chi}\right\rVert_{L^2}^2 + \left\lVert \kappa_h * \hat{\chi} - u^l \right\rVert_{L^2}^2 ,
\end{equation}
with regularization parameters $\lambda_l >0$.

Of the problem in \eqref{eq:second:cont} we consider the same discretization as in the first stage, i.e., for each angle $\theta_l$ we reconstruct a discretized version $(\chi_l )_{j,k}\in\mathbb{R}^{N_y\times N_z}$ of the X-Ray projections. To this aim we consider a discretization on the $N_y\times N_z$ grid of the $\kappa_h$ (with the midpoint rule) and denote with $K_h$ the matrix representing the convolution with $\kappa_h$. If moreover, $D$ is the matrix representing the discretization of the gradient with forward differences, we obtain the discrete version of the problem
\begin{equation}\label{eq:second:disc}
	\chi_l = \arg\min_{\hat{\chi}} \lambda_l \left\lVert D\hat{\chi}\right\rVert_2^2 + \left\lVert K_h\hat{\chi} - u^l\right\rVert_2^2 ,
\end{equation}
which can be solved by considering the corresponding Euler-Lagrange equations
\begin{equation}\label{eq:secon:EL}
	-\lambda_l L\chi_l + K_h^T \left ( K_h\chi_l-u^l\right )=0
\end{equation}
where $-L = D^T D$ is the negative Laplacian with zero Dirichlet boundary conditions.

\paragraph{\textbf{Third Stage}: Radon Inversion of the X-Ray Sinograms}

Given the reconstruction of the X-Ray projections $X_{\theta_l}[\rho ] = \chi_l$ in the second stage, we finally reconstruct the distribution. The idea is to proceed slice-wise along the $z$-axis using FBP as in Theorem \ref{thm:second:formula}. In particular, we consider the planes $\mathbf{M}_k \coloneq \lbrace (x,y,z)\in\mathbb{R}^3 \colon z=z_k\rbrace$ with $z_k$ as in \eqref{eq:disc} and we obtain the $k$-th slice of $\rho_{\mathrm{rec}}$ using a discretized implementation of \eqref{eq:second:formula}, i.e., 
\begin{equation}\label{eq:fb:disc:slice}
	\rho_{\mathrm{rec}} (\bm{x} )\bigl\lvert_{\mathbf{M}_k} = \frac{\pi}{\lvert\Theta\rvert}\sum_{l=0}^{\lvert\Theta\rvert -1}\mathcal{F}^{-1}\left ( \mathcal{F}\biggl (\chi_l\bigl\lvert_{\mathbf{M}_z}\biggr ) (\omega )\cdot \lvert\omega\rvert\cdot w (\lvert\omega\rvert )\right ) (\bm{x})
\end{equation}
for some window function $w$. A pseudocode of the complete three-stage algorithm can be found in Algorithm \ref{alg:3stage}.

\begin{rmk}\label{rmk:resolution}
	From Theorem \ref{thm:second:formula} we observe that the method inherits the same sampling-dependent resolution from the theory of the Radon transform~\cite[Chapter 3]{NattererCTbook}; in particular, if we want to resolve details of size $2\pi /b$ for some $b>0$, then one requires a number of angles $\lvert\Theta\rvert > b$ and optimally a reconstruction with $N_x = N_y = 2q+1$ with $q = \lvert\Theta\rvert/\pi$.
\end{rmk}

\begin{algorithm}[t]
	\caption{Proposed three-stage model-based reconstruction algorithm for 3D Field-Free Line Magnetic Particle Imaging\label{alg:3stage}}
	\textbf{Input}: A set of scanning angles $\Theta$ and a domain 3D $\Omega$; for each $\theta\in\Theta$, time independent samples $\bm{s}_\theta^{\theta}=\bm{s}^\theta (t_m )$, functions $\Lambda_1$,$\Lambda_2$, gradient strength $G$ and amplitudes $A_1$, $A_2$ and a regularization parameters $\lambda_l >0$ for $l=0,\dots ,\lvert\Theta\rvert -1$. \\
	\textbf{Output}: Reconstructed 3D particle density $\rho_{\mathrm{rec}}$.\\
	\begin{algorithmic} 
		\State{Initialize $\rho$ given the grid on $\Omega$;} 
		\State{Initialize scanning positions $\bm{r}_m \gets \left (\frac{A_1\Lambda_1 (t_m )}{G}\, ,\frac{A_2\Lambda_2 (t_m )}{G}\right )$ and $\bm{v}_m =\frac{d}{dt}\bm{r}_m$;}
		\For{$\theta$ in $\Theta$}
		\State{$l\gets$ index of $\theta$;}
		\State{Transform $\bm{s}_m^\theta\gets \tilde{\bm{s}}^\theta (t_m )$ as in \eqref{eq:tr:main};}
		\State{\textbackslash\textbackslash \tt{First stage: obtain the traces of the X-Ray projections.}}
		\State{Obtain MPI Core Operator $A^l\gets$ using \eqref{eq:fisrt:stage}.}
		\State{Get the traces $u^l\gets\trace (A^l )$;}
		\State{\textbackslash\textbackslash \tt{Second stage: obtain the X-Ray projections solving Eq.~\eqref{eq:second:disc}}.}
		\State{Update $\chi_l\gets\mathrm{ConjugatedGradient}(	-\lambda_l L + K_h^T  K_h\, , K_h^T u^l)$;}            
		\EndFor
		\State{\textbackslash\textbackslash \tt{Third stage: reconstruct $\rho$ in a slice-wise fashion}.}
		\For{$k\gets 0$ \textbf{to} $N_z -1$}
		\State{Update $\rho_{\mathrm{rec}} (\bm{x} )\bigl\lvert_{\mathbf{M}_k}$ using \eqref{eq:fb:disc:slice};}
		\EndFor
		\Return{The reconstruction $\rho_{\mathrm{rec}}$.}
	\end{algorithmic}
\end{algorithm}

\section{Discussion on the Generality and Feasibility of the Proposed Set-Up and the Reconstruction Formulae}\label{sec:discussion}

In this section we discuss in which way the proposed reconstruction formulae and ensuing reconstruction algorithm offer flexibility in the scanning setup. Moreover, we provide a comparison with the iMPI scanner~\cite{vogel2023impi} based on the TWMPI approach~\cite{vogel2014TWMPI} which is a good candidate for an existing scanner that could be potentially modified to realize and test the suggested set-up. 
We start with the three observations
\begin{itemize}
	\item[(O1)] the thesis in Theorem \ref{thm:main} is independent of the choice of $\Lambda_1$, $\Lambda_2$ in \eqref{eq:dynamic:field}. It is possible to consider any trajectory $(r_1 (t), r_2 (t))$ in the plane and have $\bm{r}(t) \coloneq \left (		0 \, , \frac{r_1 (t)}{G} \, , \frac{ r_2 (t)}{G} \right )^T$ in \eqref{eq:lissajous:plane}.
	\item[(O2)] In the first stage of the reconstruction algorithm the data points are set of triples of the form $\lbrace (\bm{s}^{\theta}_m , \bm{r}_m , \bm{v}_m )\colon m=0,\dots ,L_l-1\rbrace$ and in particular, the reconstruction problem in \eqref{eq:fisrt:stage} is also independent from the particular choice of $\bm{r}(t)$, provided that is covers the region of interest sufficiently.
	\item[(O3)] Given a set of X-Ray projections $\chi_l$ on a set of planes whose intersection is the z-axis, then we can apply the Second Reconstruction Formula in \eqref{eq:second:formula}.
\end{itemize}
From these observations, it follows that our reconstruction strategy is applicable for any FFL scanning set-up that, given a set of angles $\Theta\subset [0,\pi )$ and planes $\Pi_{\theta}=\langle\bm{e}_\theta\rangle^\perp$ with $\theta\in\Theta$, obeys to the properties
\begin{itemize}
	\item[(P1)] for each $\theta$, the scan is performed with the FFL parallel to $\bm{e}_\theta$;
	\item[(P2)] the trace of the trajectory $\mathrm{im}(\bm{r})$ forms an $\epsilon$-net (finite in the discretized case), i.e., \\ $\Omega\subset\bigcup_{x\in\mathrm{im}(\bm{r})} B_{\epsilon}(x)$ where $B_{\epsilon}(x)$ are open balls of radius $\epsilon$ centered in x;
	\item[(P3)] the planes form a coaxial set along the z-axis, i.e., $\cap_{\theta\in\Theta}\Pi_\theta = \langle\bm{e}_z\rangle$.
\end{itemize}
This set of properties is sufficient in the sense that, (P1) and (P2) guarantee that the first two stages are applicable for every angle $\theta$ due to (O1) and (O2). On the other hand, (P3) allows us to use the third reconstruction formula (due to (O3)). We observe that the proposed setup with Lissajous curves satisfies (P1)-(P2), but is not the only possibility.

We now briefly introduce the concept behind the TWMPI~\cite{greiner2022twmpiFFL,vogel2014TWMPI} and in particular the iMPI scanner~\cite{vogel2023impi}. In the iMPI scanner two overlapping saddle-coil pairs in Helmholtz configuration (CH1 and CH2) with a phase shift of $90^\circ$ and the same frequencies $f_1$ are used to generate an FFL that travels along the symmetry axis (y-axis). Additionally, a pair of solenoids (CH3) is able to steer the FFL along some trajectory in the xz-plane. In the general case scenario, the frequencies in the channel are chosen to be as in Tab. \ref{tab:sample-table}.
\begin{table}[H]
	\centering
	\begin{tabular}{c c c c}
		Channel  &  CH1 & CH2 & CH3 \\ \hline
		Excitation& $A\sin(2\pi f_1 )$ & $A\sin(2\pi f_1 +\pi /2 )$ & $B\sin (2\pi f_2 )$\\
	\end{tabular}
	\caption{Channels and frequencies for the iMPI scanner \cite{vogel2023impi}.}
	\label{tab:sample-table}
\end{table}
\noindent
In particular, for the case of the iMPI scanner \cite{vogel2023impi}, $f_1 = 60$ Hz and $f_2 = 2480$ Hz, with the results that the FFL is driven along a sinusoidal path perpendicularly to the xz-plane. This is the case because $f_1 \ll f_2$, however, it is enough to consider $f_1 = f_2$ with $f_1 / f_2 = p/(p+1)$ for some $p\in\mathbb{N}$ to obtain a closed Lissajous curve. It is however not necessary, because the reconstruction procedure is independent from the chosen scan trajectory. Finally, the iMPI scanner can be rotated along an axis to obtain multiple projections and reconstruct a 3D object as proposed in this paper.

\section{Numerical Experiment}\label{sec:experiments}

In this section we describe a numerical experiment to demonstrate the potential of the reconstruction algorithm in Section \ref{sec:rec:algorithm} for a simulated 3D FFL scan.

\paragraph{\color{black}{Experimental Setup and Simulation of the Scans}}

In the simulated experiment, we consider a scanner that produces magnetic fields as in \eqref{eq:magfield} with parameters inspired by \cite{thieben20243dmpi}. More specifically, the scanner in \cite{thieben20243dmpi} is able to produce 2D FFP Lissajous scans with components $r_i (t) = A_i \sin (2\pi f_i t)$ for $i=1,2$ with $f_1 = \frac{125}{4864}$\si{\mega\hertz} and $f_1 = \frac{125}{4800}$\si{\mega\hertz}, yielding a frequency ratio of $f_1 / f_2 = 75/76$ and a total period of $T_\mathrm{loop} = 2.9184$\si{\milli\second} for a full scan. The sampling frequency is $f_\mathrm{sample} = \frac{64}{125}$\si{\mega\hertz} which yields a total of $L = T_\mathrm{loop}\cdot f_\mathrm{sample} = 5700$ sample points. We imagine a scanner that produces 2D Lissajous trajectories similarly to \cite{thieben20243dmpi} (but with a FFL extending in the third direction); in particular, we consider a gradient $\bm{G} = \diag (-0.12 ,0.12, 0 ) \frac{T}{\mu_0 m}$, and excitation of about 5\si{\milli\tesla}, yielding amplitudes $A_i = -33.\bar{3}\si{\milli\meter}$ and consequently a 2D FoV covered by the Lissajous curve is about $67\times 67 \si{\milli\meter}^2$. Concerning the particle distribution, we have simulated a phantom that resembles a stenotic blood vessel (cf. Fig.\ref{subfig:rec:gt}) on a $500\times 500 \times 500$ grid and filled with \emph{perimag}\textregistered plain nanoparticles (micromod Partikeltechnologie, Germany) with hydrodynamic diameter of $130\si{\nano\meter}$. Assuming the particles to be at the average body temperature of $310\si{\kelvin}$ during the scan, we obtain $H_\mathrm{sat}\approx 23.24\si{\ampere / \meter}$. With the stated size of the FoV and strength of magnetic field, the resolution parameter is $h\approx 0.00365 \si{\ampere / \meter}$. For the simulation of the signal $\bm{s}^{\theta}$ we have used the relationship in \eqref{eq:main} with a trivial coil sensitivity pattern  $P=\mathrm{Id}_3$ the identity matrix.
We consider scans with $100$ equidistant angles $\theta_l = l\cdot \frac{\pi}{100}$ for $l = 0,\dots , 99$ and for each angle $\theta_l$ we have considered $L = 5700$ time samples $t_m = m\cdot \frac{T}{L}$ for $m=0, \dots ,L-1$, assuming $T=T_\mathrm{loop}$ for each angle. It follows that if $T_\mathrm{rot.}$ is the time needed for an angle rotation, the total duration of a scan is $L\cdot T_\mathrm{loop}+(L-1)\cdot T_\mathrm{rot.}$.

To account for noise contamination, we have corrupted the simulated signals with $2\%$ additive Gaussian noise, i.e., for each angle $\theta$ we considered a corrupted signal
\begin{equation}\label{eq:noisy:data}
	\hat{\bm{s}}^{\theta} (t_m ) = \bm{s}^{\theta}(t_m )+\varepsilon_\theta N_m^{\theta} \quad\theta \in\left\lbrace l\cdot \frac{\pi}{100}\colon l = 0,\dots , 99\right\rbrace .
\end{equation}
Here $N_m^{\theta}$ are i.i.d.~random variables extracted from a normal distribution with zero mean and standard deviation of one. The noise amplitude in this case is 
\begin{equation}
	\varepsilon_\theta = 0.02 \cdot\max_m \lbrace\lvert \bm{s}^{\theta}(t_m )\rvert\rbrace .
\end{equation}
The signal in \eqref{eq:noisy:data} serves as input of the reconstruction algorithm in Section \ref{sec:rec:algorithm}.

\paragraph{\color{black}{Reconstruction Parameters}}

To avoid reconstructing on the same grid used for the simulation, the reconstruction has been performed on a $50\times 50\times 50$ grid. In particular, for every angle $\theta_l$ we have reconstructed the trace of the MPI Core Operator $u_l$ variationally as described in Section \ref{sec:rec:algorithm} on a $50\times 50$ grid. In particular, the Euler-Lagrange equations have been solved using the CG method with a parameter $\mu = 18$ obtained with a grid search optimizing with respect to the Peak-Signal-to-Noise (PSNR) metric. Each of these traces $u_l$ have been deconvolved by solving with the CG method the Euler-Lagrange equations associated to the deconvolution problem in \eqref{eq:secon:EL}, with parameter $\lambda = \lambda_l = 5\times 10^{-5}$ equal for all angles and fine tuned and obtained with a grid search optimizing with respect to PSNR. The maximum number of iterations for the CG algorithms in the first and second stages have been set to $10000$, if a relative tolerance of $10^{-6}$ has not been reached. In particular, all reconstruction converge below tolerance before reaching the $10000$ iterations. In order to introduce an additional form of unpredictable noise, we consider for the deconvolution $\tilde{h}=0.004$, introducing a $12\%$ relative error on $h$. The set of reconstructed X-Ray projections $\chi_l$ serves as the input for the third and last stage, where we employ the FBP formula in \eqref{eq:secon:EL} on each slice perpendicular to the $z$-axis, with the Ramp-filter $w(x)=x$ if $x\leq 0.5$ and $w(x)=0$ otherwise.

The algorithm described in the paper and the data simulation have been implemented in Python 3.9, using the packages Numpy, SciPy, Scikit-image and PyTorch. The numerical experiments were performed on a desktop computer with 13th Gen Intel(R) Core(TM) i9-13900KS, 128 GB of RAM, an NVIDIA RTX A6000 GPU and Windows 11 Pro.

Finally, as a study on the effect of different trajectories and comparison with the TWMPI approach, we have simulated the scan produced with the frequencies of the iMPi scanner $f_1 = 60 \si{\hertz}$ and $f_2 = 2480\si{\hertz}$ with ratio of $f_1 / f_2 = 60/2480 = 3/124$ and keeping all other parameters as described. For this scan, the parameters obtained with the same grid search are $\mu = 2$ and $\lambda = 3\cdot 10^{-4}$.

\def\imratio{0.32}
\begin{figure}[h]
	\centering
	\begin{subfigure}[t]{\imratio\linewidth}
		\includegraphics[width=\linewidth]{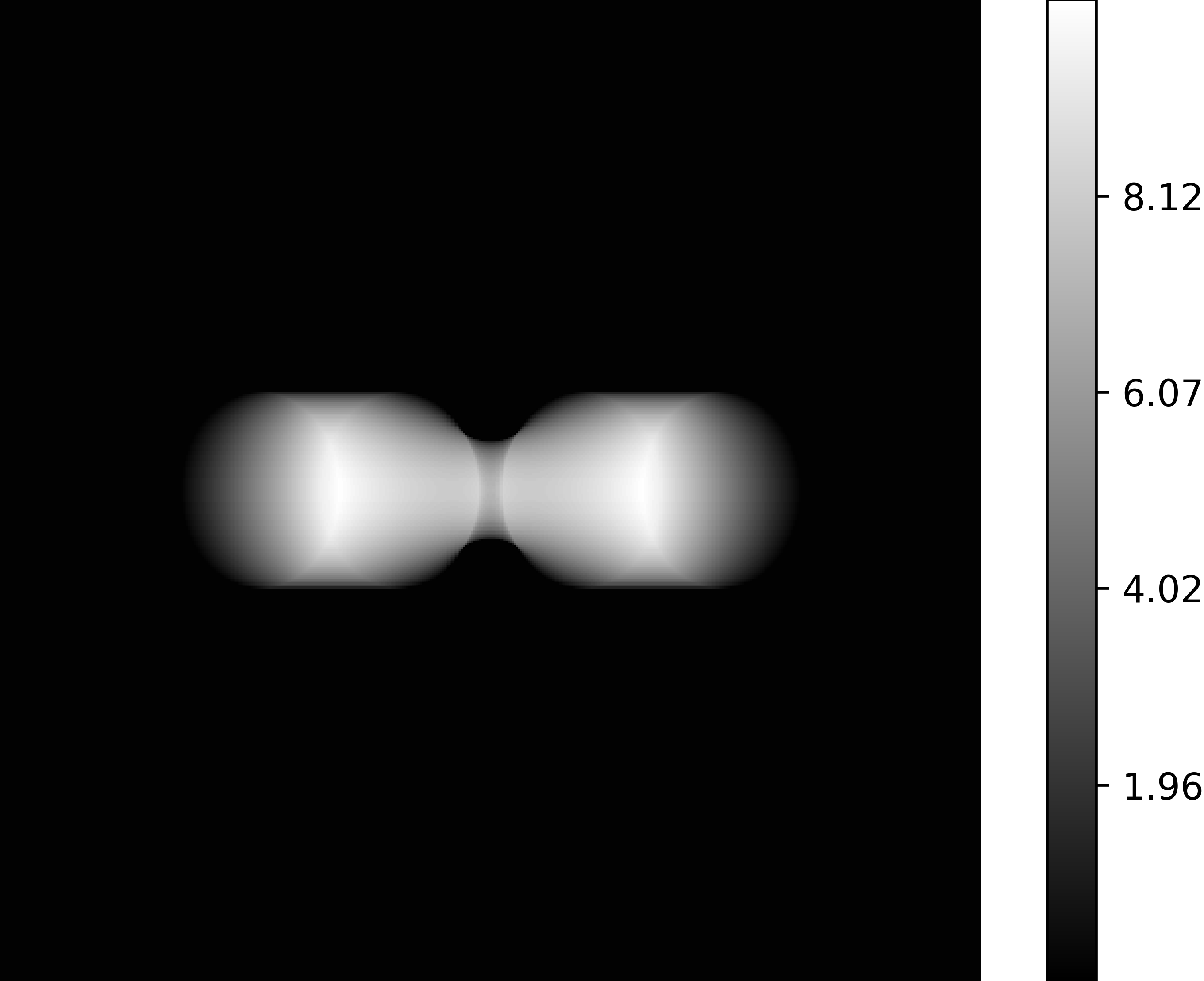}
		\caption{\centering\scriptsize X-Ray Projection $144^\circ$}
		\label{subfig:2stages:gt}
	\end{subfigure}
	\hfil
	\begin{subfigure}[t]{\imratio\linewidth}
		\includegraphics[width=\linewidth]{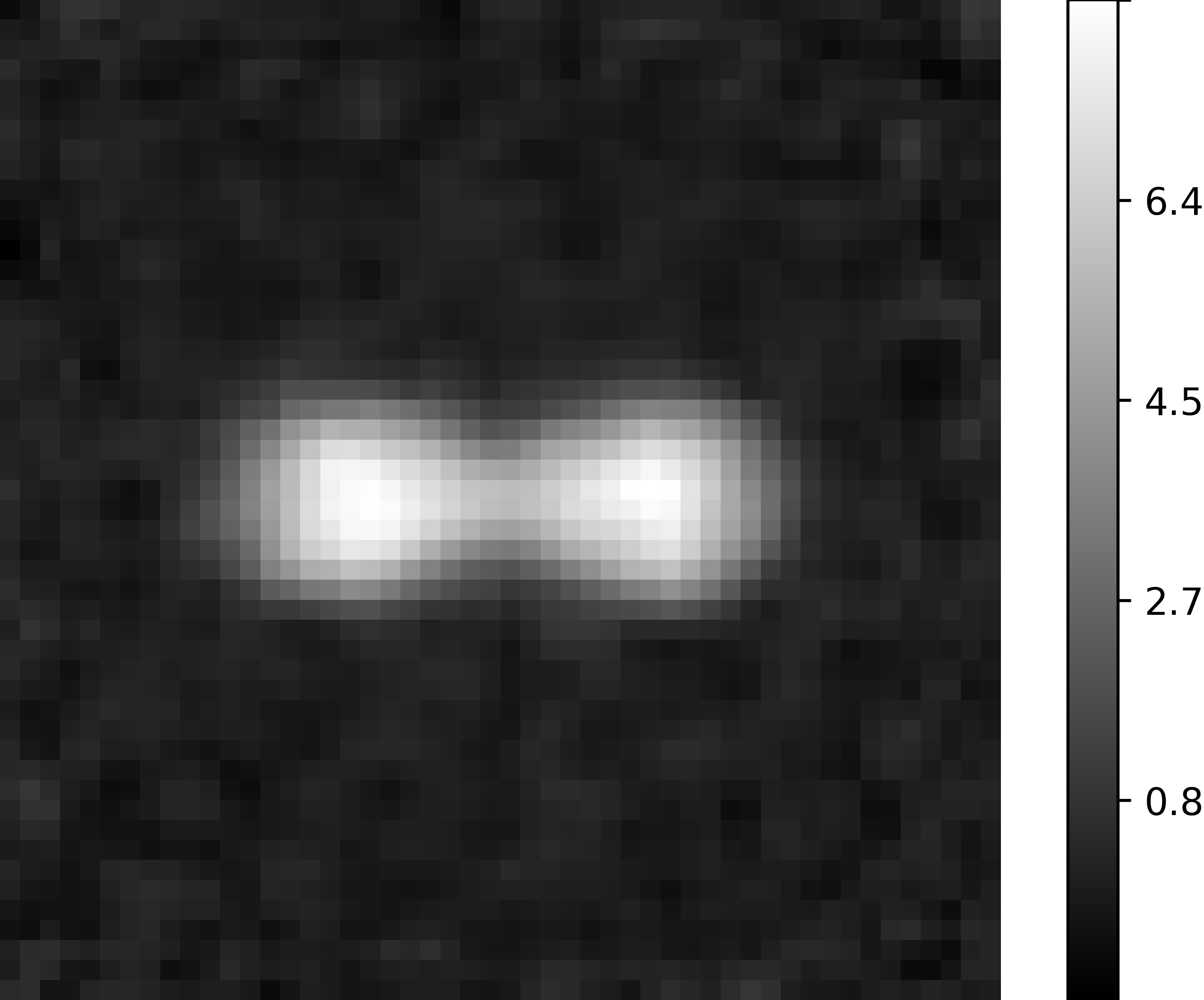}
		\caption{\centering\scriptsize Proposed reconstruction of $ X_{\bm{e}_\theta} [\rho ]$. }
		\label{subfig:2stages:proposed}
	\end{subfigure}
	\hfil
	\begin{subfigure}[t]{\imratio\linewidth}
		\includegraphics[width=\linewidth]{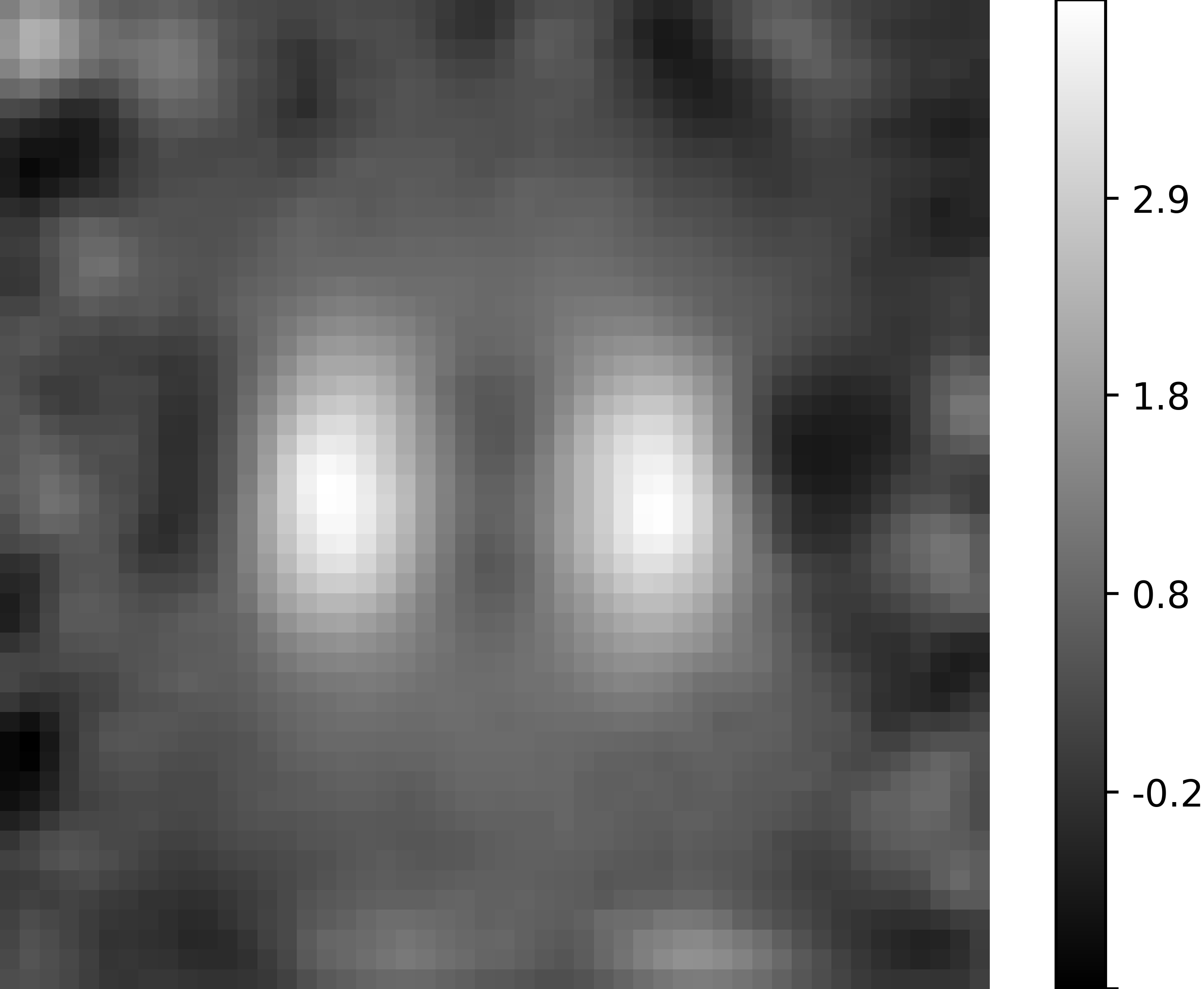}
		\caption{\centering\scriptsize iMPI reconstruction of $ X_{\bm{e}_\theta} [\rho ]$. }
		\label{subfig:2stages:iMPI}
	\end{subfigure}
	\par
	\begin{subfigure}[t]{0.49\linewidth}
		\centering
		\includegraphics[width=\linewidth]{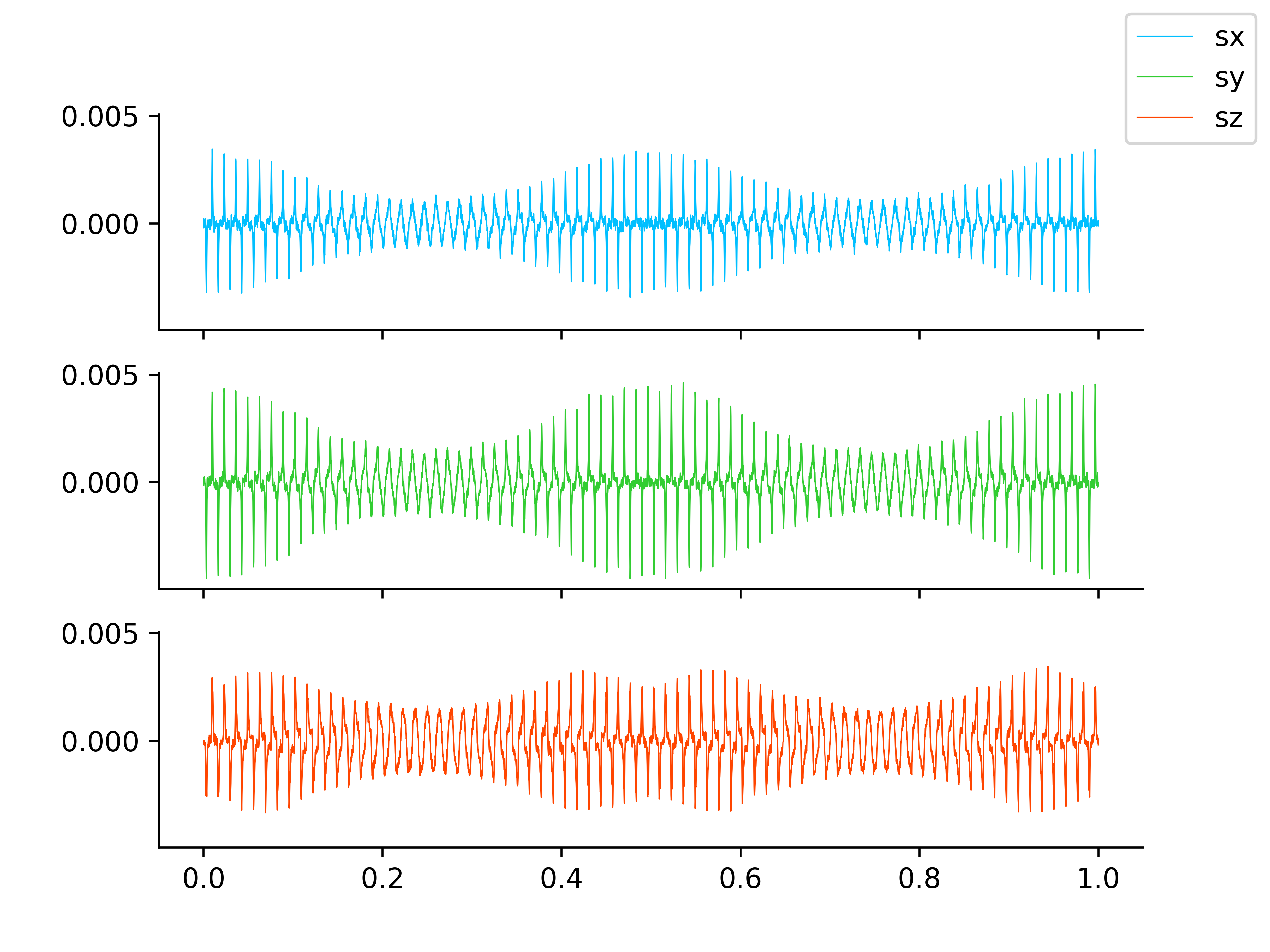}
		\caption{\centering\scriptsize Scan with $f_1 / f_2 = 76/75$.}
		\label{subfig:2stages:plot}
	\end{subfigure}
	\hfil
	\begin{subfigure}[t]{0.49\linewidth}
		\centering
		\includegraphics[width=\linewidth]{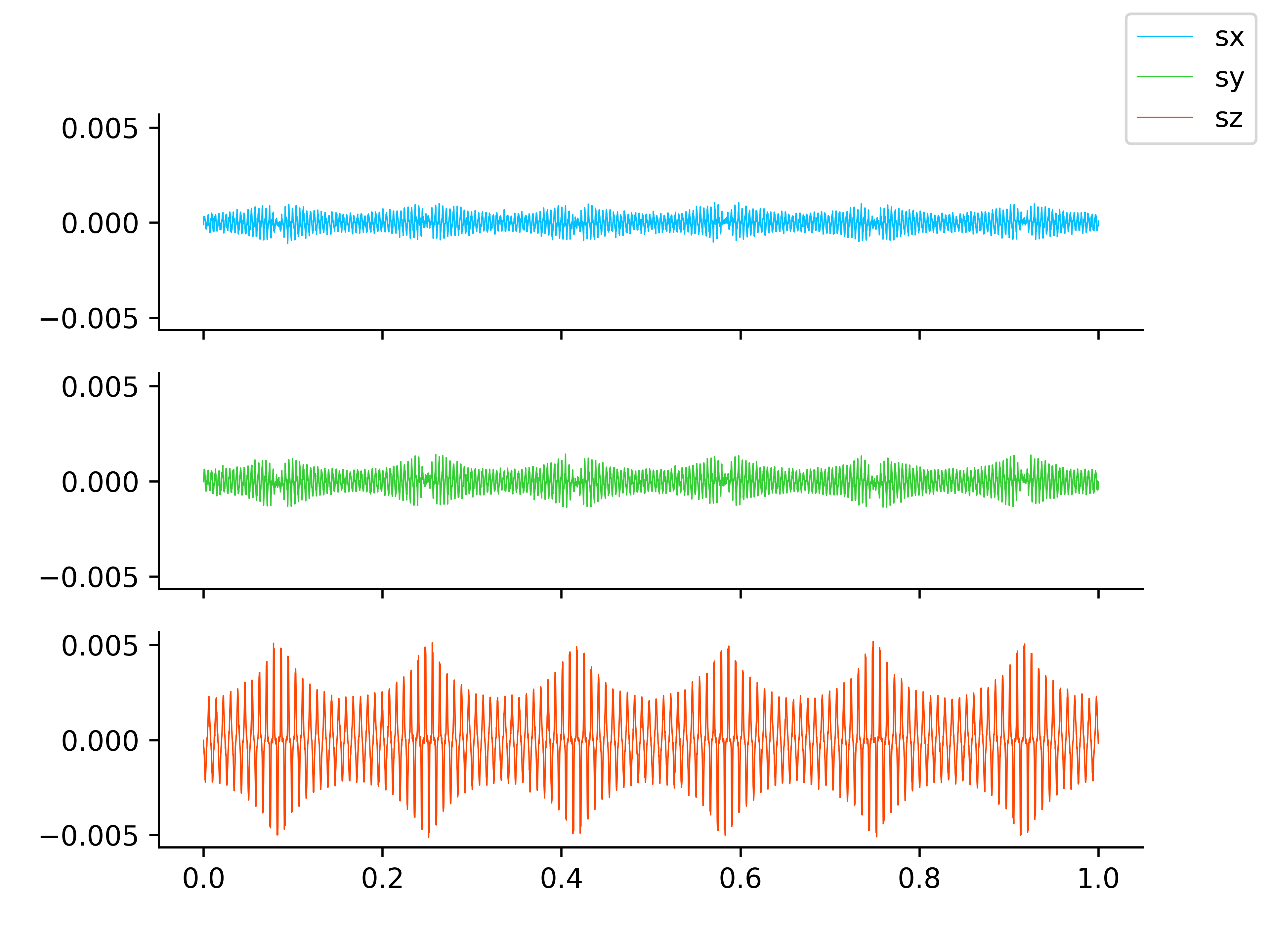}
		\caption{\centering\scriptsize  Scan with $f_1 / f_2 = 3/124$.}
		\label{subfig:2stages:plot2}
	\end{subfigure}
	\caption{(\textbf{a}) Ground truth of the X-Ray Projection for $\theta = 144^\circ$. \textbf{(b)} Reconstructed projection using the proposed scanning setup using the scan data in (d). \textbf{(d)} Reconstructed projection with the iMPI setup using the scan data in (e). \textbf{(d)-(e)} Comparison of the the components along the x, y and z (from top to bottom) of the signals $\bm{s}^{\theta}(t_m )$ obtained during the scan with $\theta =144^\circ $ with both the proposed and the iMPi methods.  }
	\label{fig:2stages}
\end{figure}

\def\imratio{0.30}
\begin{figure}[H]
	\centering
	\begin{subfigure}[t]{\imratio\linewidth}
		\includegraphics[width=\linewidth]{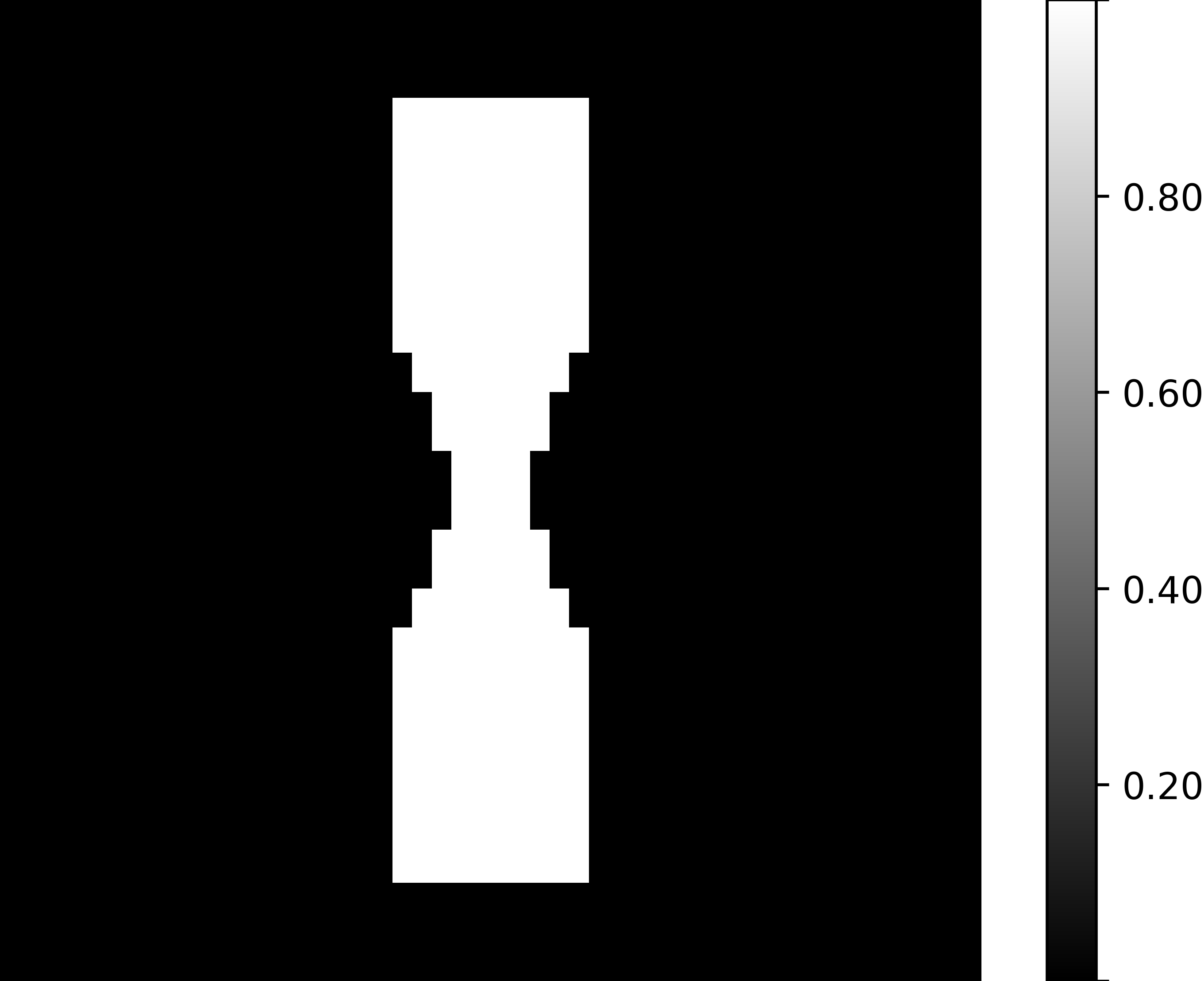}
		\caption{\centering\scriptsize $\rho_{\mathrm{GT}}$ slice}
		\label{subfig:rec:a}
	\end{subfigure}
	\begin{subfigure}[t]{\imratio\linewidth}
		\includegraphics[width=\linewidth]{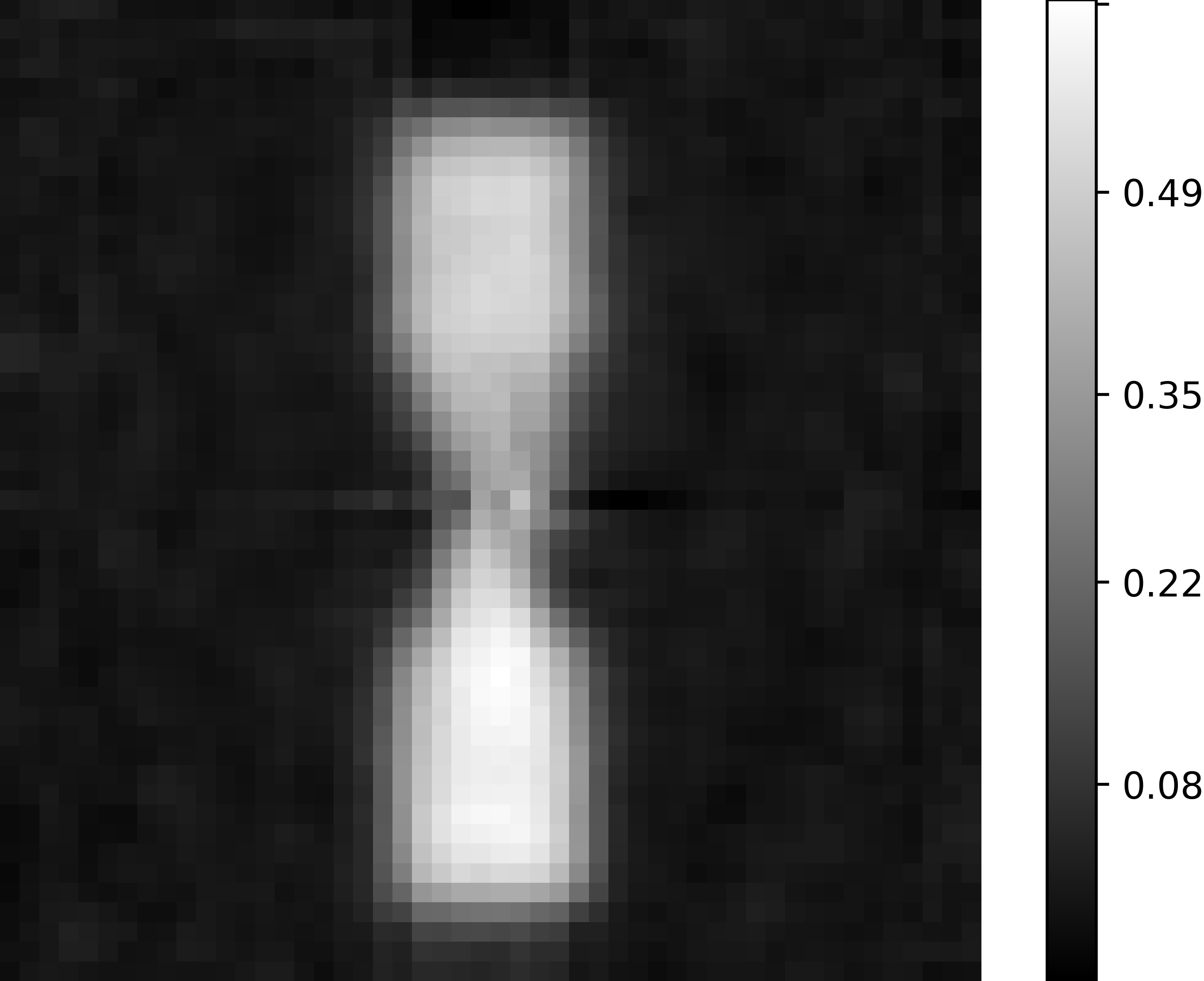}
		\caption{\centering\scriptsize $\rho_{\mathrm{rec}}$  with the proposed setup}
		\label{subfig:rec:b}
	\end{subfigure}
	\begin{subfigure}[t]{\imratio\linewidth}
		\includegraphics[width=\linewidth]{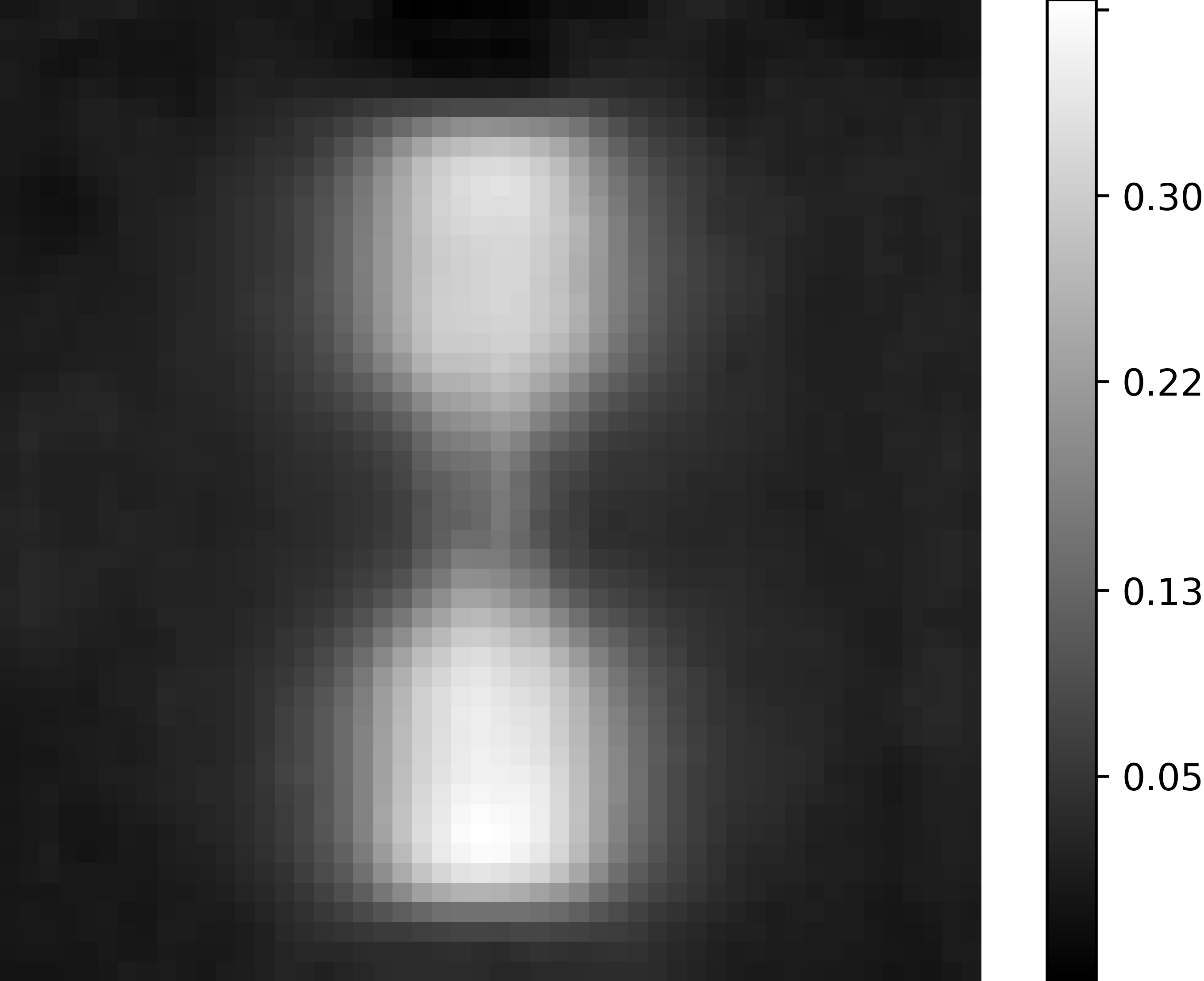}
		\caption{\centering\scriptsize $\rho_{\mathrm{rec}}$ 
			with the iMPI-type setup}
		\label{subfig:rec:c}
	\end{subfigure}
	\centering
	\begin{subfigure}[t]{0.32\linewidth}
		\includegraphics[width=\linewidth]{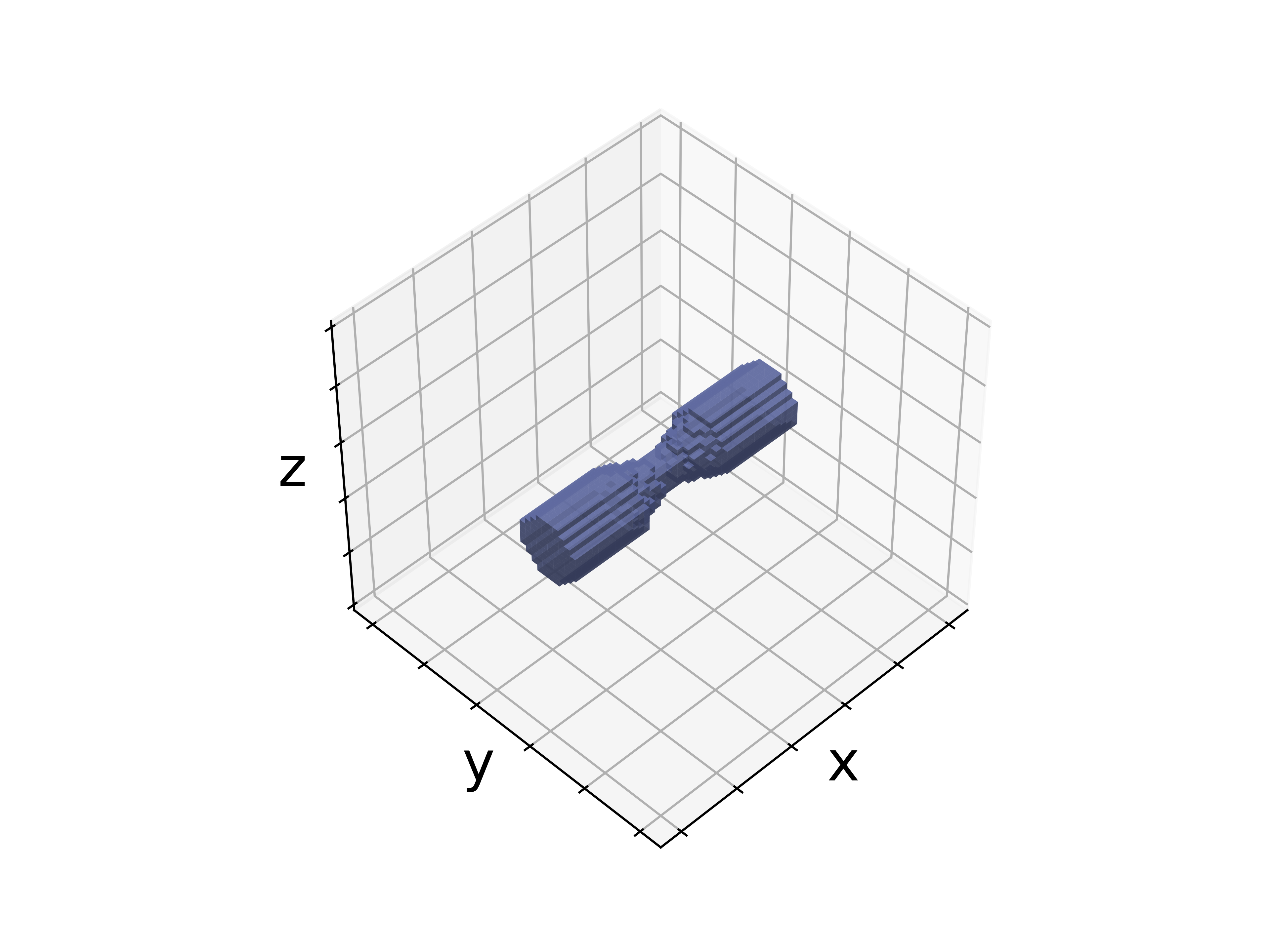}
		\caption{\centering\scriptsize $\rho_{\mathrm{GT}}$}
		\label{subfig:rec:gt}
	\end{subfigure}
	\begin{subfigure}[t]{0.32\linewidth}
		\includegraphics[width=\linewidth]{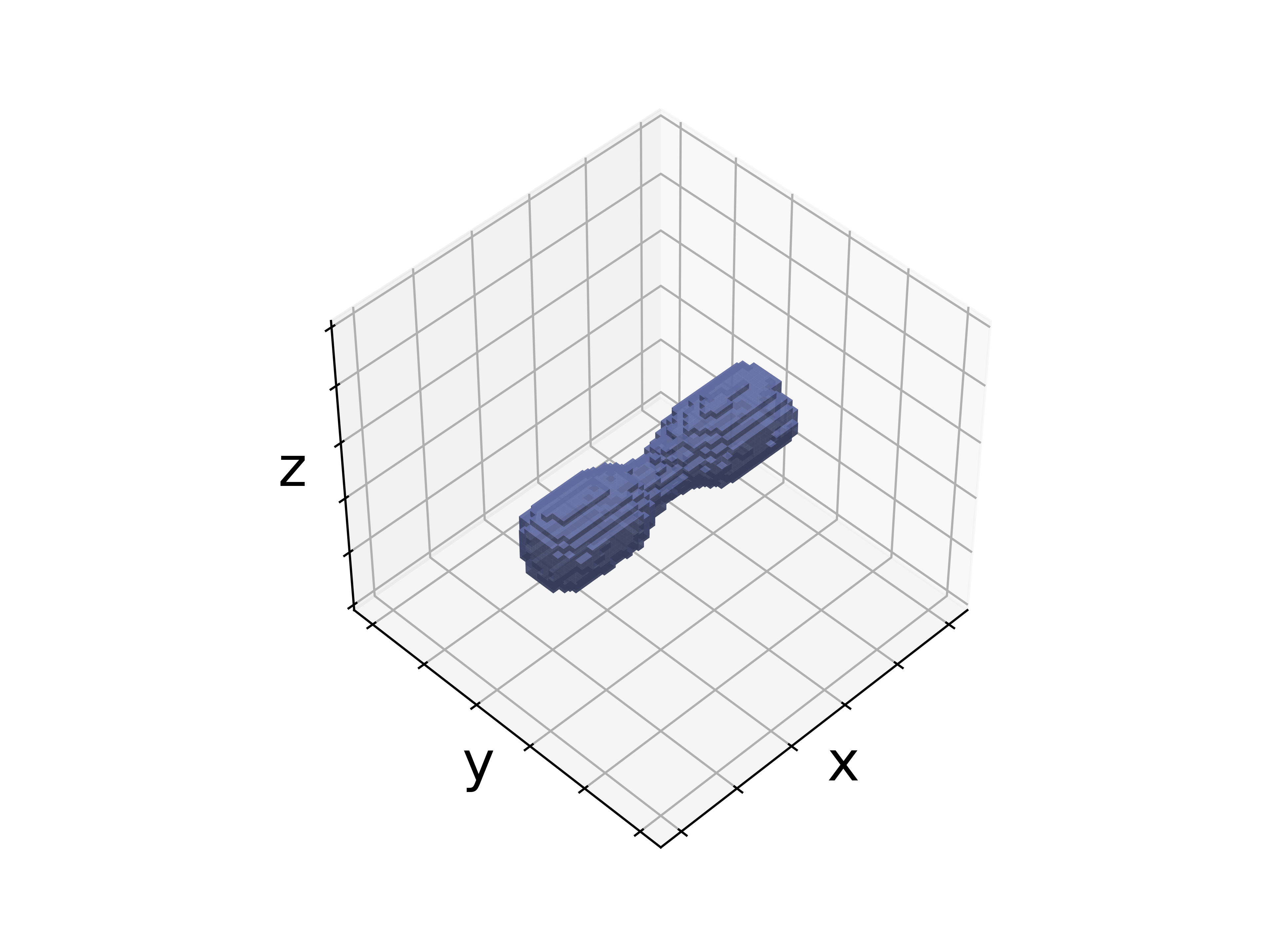}
		\caption{\centering\scriptsize $\rho_{\mathrm{rec}}$ with proposed setup}
		\label{subfig:rec:rec}
	\end{subfigure}
	\begin{subfigure}[t]{0.32\linewidth}
		\includegraphics[width=\linewidth]{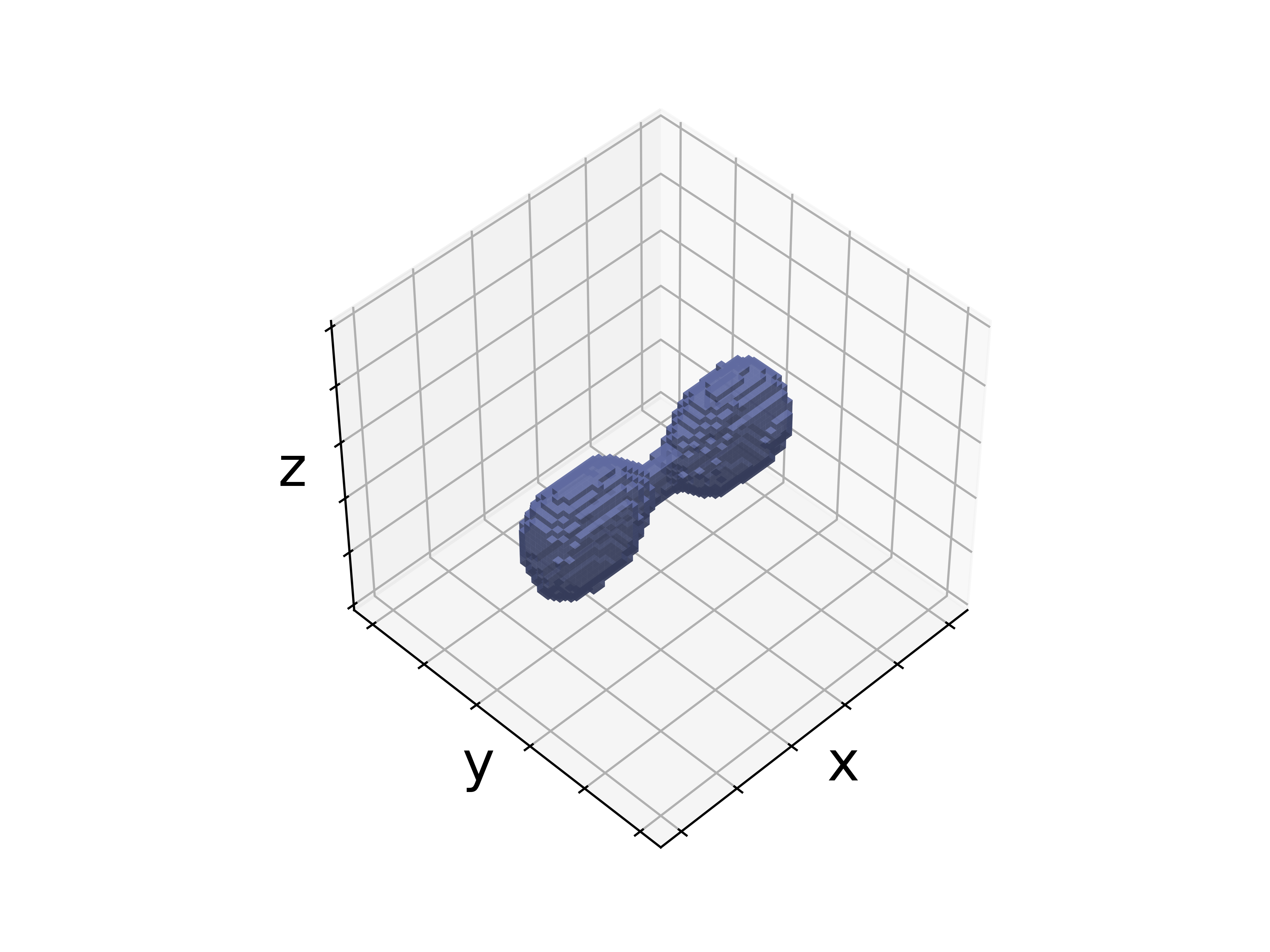}
		\caption{\centering\scriptsize $\rho_{\mathrm{rec}}$ with iMPI-type scan}
		\label{subfig:rec:rec_iMPI}
	\end{subfigure}
	\caption{
		Reconstruction of a stenotic vessel-shaped phantom with the 3D FFL scan proposed and the iMPI-inspired scan. \textbf{(a)} A selected slice of the ground truth along the y-axis; \textbf{(b)} the reconstruction of the selected slice with our algorithm and the proposed scanning setup; \textbf{(c)} the reconstruction of the slice with our algorithm and the iMPI scanning setup; \textbf{(d)-(f)} boolean plots of the ground truth $\rho_\mathrm{GT}$ and the final 3D reconstructions after Filtered Backprojection. }
	\label{fig:rec}
\end{figure}

\paragraph{Discussion of the Results} In Fig.\ref{fig:2stages} we can see an example of the reconstructions of the X-Ray projections with the the algorithm presented in this paper: a reconstruction from the signal in Fig.\ref{subfig:2stages:plot} and Fig.\ref{subfig:2stages:plot2}, obtained from the scan performed at angle $\theta = 144^\circ$. In Fig.\ref{fig:rec} we can see the final reconstructions after the third stage of the algorithm (the Radon inversion). The reconstructed distribution is plotted in boolean values in Fig.\ref{subfig:rec:rec} and Fig.\ref{subfig:rec:rec_iMPI}, where we put all values below $25\%$ of $\max\lbrace\lvert\rho_{\mathrm{rec}}\rvert\rbrace$ to 0 and all the others to 1. For a fairer comparison between the $\rho_{\mathrm{GT}}$ and $\rho_{\mathrm{rec}}$ we have also plotted an example slice. In particular, Fig.\ref{subfig:rec:a}, \ref{subfig:rec:b} and \ref{subfig:rec:c} are the ground truth, the reconstruction with $f_1 /f_2 = 76/75$ and with $f_1 / f_2 = 3/124$, respectively. Comparing the ground truth and the reconstructions it is evident that the main algorithm is capable of reconstructing shapes, but the edges of the reconstruction are smoothed out. This smoothing out of the reconstruction can be already seen in the reconstructed X-Ray projections, for example in Fig.\ref{subfig:2stages:iMPI}. This is due to the fact that indeed the regularizer (see Fig.\ref{eq:second:cont}) employed in the second stage enforces smoothness of the reconstruction and could be mitigated employing regularizers that are more suitable for the reconstruction of edges (more details in Section \ref{sec:conclusions}). This experiment demonstrates that a model-based reconstruction of 3D distributions from 3D FFL scans using our algorithm and methodology is promising.
\def\imratio{0.32}
\begin{figure}[h]
	\centering
	\begin{subfigure}[t]{\imratio\linewidth}
		\includegraphics[width=\linewidth]{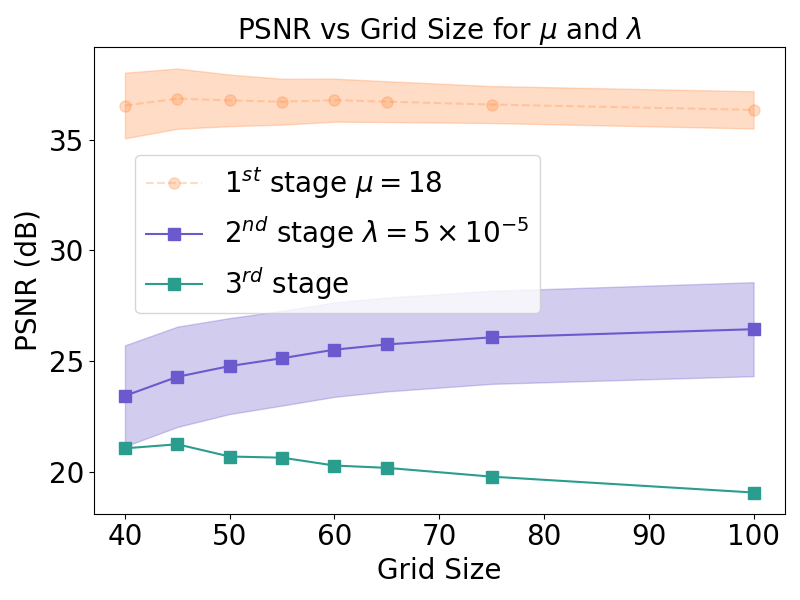}
		\caption{\centering\scriptsize Grid Size vs. PSNR}
		\label{subfig:metric:grid}
	\end{subfigure}
	\begin{subfigure}[t]{\imratio\linewidth}
		\includegraphics[width=\linewidth]{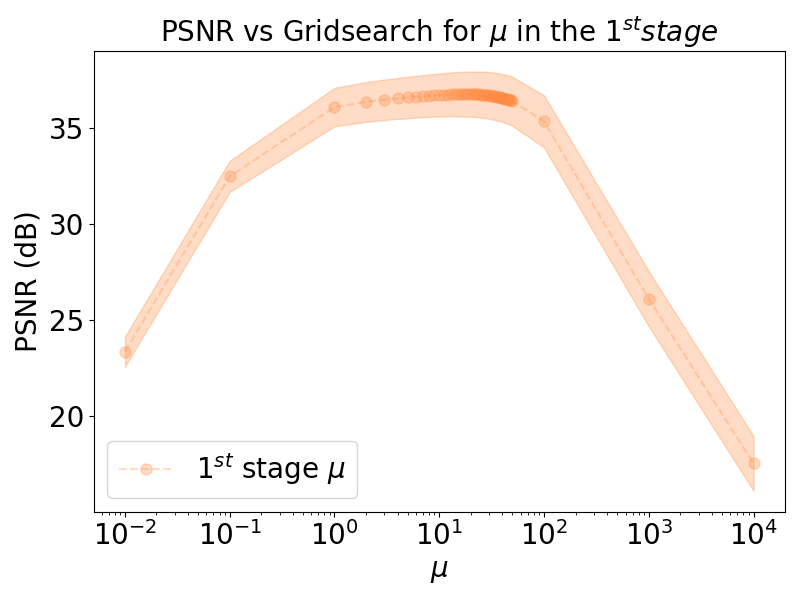}
		\caption{\centering\scriptsize $1^{st}$ stage $\mu$ vs PSNR}
		\label{subfig:metric:first}
	\end{subfigure}
	\begin{subfigure}[t]{\imratio\linewidth}
		\includegraphics[width=\linewidth]{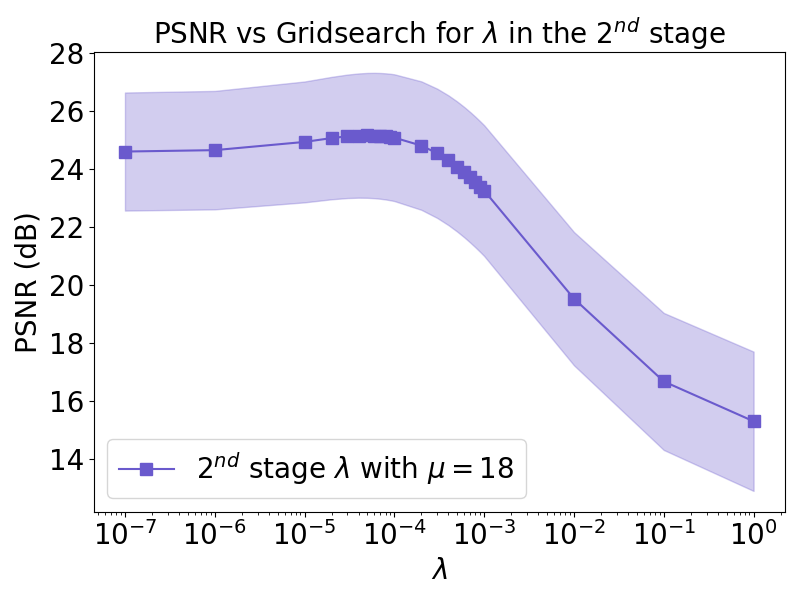}
		\caption{\centering\scriptsize $2^{nd}$ stage $\lambda$ vs PSNR}
		\label{subfig:metric:second}
	\end{subfigure}
	\caption{Resolution and stability of the reconstruction algorithm with respect to the grid size, $\mu$ and $\lambda$ where the thick lines represent the averages and the shaded areas show one standard deviation over the $100$ angles. \textbf{(a)} PSNR obtained by keeping fixed the parameters $\mu = 18$, $\lambda = 5\cdot 10^{-5}$ and varying the reconstruction grid size $N_x = N_y$ from 40 to 100. \textbf{(b)} PSNR computed on the values of the grid search for $\mu$ in the first stage for the $50\times 50$ reconstruction grid; \textbf{(c)} PSNR computed performing the grid search for $\lambda$ with $\mu = 18$ fixed.}
	\label{fig:metrics}
\end{figure}
As an analysis on the pixel resolution, we have performed reconstructions for grid sizes ranging from $40\times 40$ to $100\times 100$ while keeping $\mu = 18$ and $\lambda = 5\cdot 10^{-5}$ fixed. The PSNR plots are displayed in Fig. \ref{subfig:metric:grid}, from which we observe that the first stage is rather stable for various choices of $\mu$, whereas the second stage even benefits from a denser grid, whereas the PSNR computed on the 3D reconstruction after Filtered Backprojection decreases; this is to be expected in view of Remark \ref{rmk:resolution}, as the reconstruction on a $N\times N$ grid, one has $q=\frac{N-1}{2}= \lvert\Theta\rvert/\pi$ for $\lvert\Theta\rvert$ angles. It follows for $\lvert\Theta\rvert=100$ that $N = \frac{2\lvert\Theta\rvert}{\pi}+1\approx 65$. Consequently, with a grid of size greater than $N_x = 65$, $100$ angles is below the optimal number of angles required. In Fig. \ref{subfig:metric:first} we plot the behavior for the PSNR metric on the grid search for the parameter $\mu$ and in Fig. \ref{subfig:metric:first} the PSNR behavior for $\lambda$ with fixed $\mu = 18$. From Fig. \ref{subfig:metric:first} and \ref{subfig:metric:second} we observe that the quality of the reconstruction measured in PSNR is rather stable under reasonable perturbation of the optimal parameters found with the grid searches.

\section{Conclusions}\label{sec:conclusions}

In this work we have provided model-based reconstruction formulae for 3D FFL MPI that relate the signal obtained in the receiving coils to the target distribution via the MPI Core Operator and the X-Ray projection. Based on those, as our second contribution we have proposed a three stage reconstruction algorithm: for each angle, the first two stages are employed to first reconstruct the traces of the MPI Core Operators involved and then, to reconstruct the X-Ray projections of the ground truth. In the third stage, the whole 3D distribution is reconstructed using the Fourier Slice Theorem. Finally, we have demonstrated the applicability of the proposed algorithm with a simulated numerical example and compared it with a simulation of a modified version of the iMPI scanner \cite{vogel2023impi}. Directions of future work include the development and testing of tools for the enhancement of each of the three stages of the algorithm presented, for example methods that encourage sparsity of the reconstruction and enforce positivity~\cite{gapyak2023multipatch} and methods that have better edge preserving properties, like TV type regularizers~\cite{rudin1992nonlinear} as well as Potts type regularizers~\cite{weinmann2015iterative,storath2015joint,weinmann2013l1potts}. Another important direction of research is the application of the algorithm to real data. It is known that the LFVs in the FFL acquisition topology are not ideal lines and rather resemble 3-dimensional banana-shaped cylinders~\cite{bringout2020new}. In order to counteract artifacts coming from the assumption of ideal FFLs, work should be also done to incorporate existing more realistic models~\cite{bringout2020new} into the stages of the algorithms here presented. Finally, the inclusion of relaxation effects ~\cite{Croft2012relaxation,li2023jiles-arthenton} in the model and the development of reconstruction formulae that take it into account are one of the main directions of research.

\section*{Acknowledgments}
We would like to acknowledge the support of the Hessian Ministry of Higher Education, Research, Science and the Arts within the Framework of the ``Programm zum Aufbau eines akademischen Mittelbaus an hessischen Hochschulen", by the German Science Fonds DFG under grant INST 168/4-1 and funding from EUt+/h-da (\url{https://h-da.de/hochschule/eut}).

\bibliographystyle{ieeetr}
\bibliography{literature}

\end{document}